\documentclass[10pt]{article}

\usepackage[colorlinks,linkcolor=blue,citecolor=blue]{hyperref}
\usepackage{amsthm}
\usepackage{multirow}
\usepackage{marginnote}
\usepackage{a4wide}
\usepackage{amssymb}
\usepackage{amsfonts}
\usepackage{amsmath}
\usepackage{mathrsfs}
\usepackage{tikz}
\usepackage{tikz-cd}
\usetikzlibrary{arrows,matrix}
\usetikzlibrary{positioning}
\usepackage{mdframed} 
\usepackage{lipsum} 
\usepackage{extarrows} 
\usepackage{color}
\usepackage{enumitem} 
\usepackage{tocloft} 
\usepackage{titlesec} 
\usepackage[T1]{fontenc}
\usepackage{titling}

\usepackage[all,cmtip]{xy}
\input xy
\xyoption{arrow} \xyoption{matrix}

\DeclareMathAlphabet{\mathpzc}{OT1}{pzc}{m}{it}
\DeclareMathAlphabet\EuScript{U}{eus}{m}{n}

\date{}

\newtheorem{proposition}{Proposition}[section]
\newtheorem{theorem}[proposition]{Theorem}
\newtheorem{lemma}[proposition]{Lemma}
\newtheorem{example}[proposition]{Example}
\theoremstyle{definition}\newtheorem{definition}[proposition]{Definition}
\theoremstyle{definition}\newtheorem{remark}[proposition]{Remark}
\newtheorem{corollary}[proposition]{Corollary}

\def\GK{{\rm  GK}\,}

\def\s{\sigma}
\def\Z{\mathbb{Z}}
\def\l{\lambda}

\def\d{\delta}

\def\CS{{\cal S}}

\def\CG{{\cal G}}
\def\CH{{\cal H}}

\def\i{{\bf i}}

\def\supp{{\rm supp }}

\def\ker{ {\rm ker } }

\def\D{ \Delta }

\def\CR{ {\cal R}}
\def\CR{ {\cal R}}

\def\CH{ {\cal H}}

\def\End{ {\rm End }}

\def\CR{ {\cal R }}

\def\ga{\mathfrak{a}}
\def\gb{\mathfrak{b}}

\def\gn{\mathfrak{n}}
\def\gm{\mathfrak{m}}
\def\gp{\mathfrak{p}}

\def\Max{{\rm Max}}

\def\N{\mathbb{N}}

\def\gl{{\rm gl}}

\def\mS{\mathbb{S}}

\def\ann{{\rm ann}}

\def\mT{\mathbb{T}}
\def\ind{{\rm ind}}

\def\mX{\mathbb{X}}

\def\Frac{{\rm Frac}}

\def\ga{\mathfrak{a}}

\def\tor{{\rm tor}}

\def \CB{\mathcal{B}}

\def\sl{\mathfrak{sl}}
\def\gl2{\mathfrak{gl}_2}

\def\mK{\Bbbk}
\def\Prim{{\rm Prim}}

\def\sC{\mathscr{C}}

\def \CO{\mathcal{O}}

\def\Aut{{\rm Aut}}

\def\SL{{\rm SL}}

\def\Spec{{\rm Spec}}

\makeatletter
\newenvironment{proof*}[1][\proofname]{\par
  \pushQED{\qed}%
  \normalfont \partopsep=\z@skip \topsep=\z@skip
  \trivlist
  \item[\hskip\labelsep
        \itshape
    #1\@addpunct{.}]\ignorespaces
}{%
  \popQED\endtrivlist\@endpefalse
}
\makeatother

\setlist[enumerate]{topsep=0pt,itemsep=-1ex,partopsep=1ex,parsep=1ex}

\newenvironment{psmallmatrix}
{\left(\begin{smallmatrix}}
	{\end{smallmatrix}\right)}

\renewenvironment{abstract}
{\par\noindent\textbf{\abstractname.}\ \ignorespaces}
{\par\medskip}

\begin{document}
\author{W.-Q. Tao} 
\title{The Heisenberg double of the quantum Euclidean group and its representations}
\maketitle
\begin{abstract}
The Heisenberg double $D_q(E_2)$ of the quantum Euclidean group $\CO_q(E_2)$ is the smash product of $\CO_q(E_2)$ with its Hopf dual $U_q(\mathfrak{e}_2)$.  For the algebra $D_q(E_2)$, explicit descriptions of its prime, primitive, and maximal spectra are obtained. All prime factors of $D_q(E_2)$ are presented as generalized Weyl algebras. As a result, we obtain that the algebra $D_q(E_2)$ has no finite-dimensional representations, and that $D_q(E_2)$ cannot have a Hopf algebra structure.  The automorphism groups of the quantum Euclidean group and its Heisenberg double are determined.  Some centralizers are explicitly described via generators and defining relations. This enables us to give a classification of simple weight modules, and the so-called $a$-weight modules, over the algebra $D_q(E_2)$. 

\


	
\end{abstract}

\section{Introduction}
In this paper, a module means a left module, $\mK$ is an algebraically closed field of characteristic zero, $\mK^*:=\mK \setminus \{0\}$ and $q\in \mK^*$ is not a root of unity.

The Euclidean group $E_2$ of orientation preserving isometries of the plane $\mathbb{R}^2$ is the semi-direct product of groups $SO_2 \ltimes \mathbb{R}^2$. 
The quantum Euclidean group $\CO_q(E_2)$ is obtained by the contraction of the quantum group $SU_q(2)$ in \cite{WoronowiczContra}. This is similar to the classical construction whereby $E_2$ is obtained as a limit of the simple group $SO_3$.  The idea to use contraction procedure in the theory of quantum groups goes back to \cite{Celeghini-Giachetti-Sorace-Tarlini}.  In the literature, $\CO_q(E_2)$ is interpreted as a quantization of the algebra of regular functions on $E_2$, see for instance \cite{Koelink} and \cite[Section 13.4]{Chari-Pressley}.   Woronowicz \cite{Woronowicz92,WoronowiczContra} investigated the quantum Euclidean group from the $C^*$-algebra level, he also constructed and investigated its Pontryagin dual \cite{Woronowicz91}. Koelink \cite{Koelink} and Korogodskii and Vaksman \cite{Korogodskii-Vaksman} studied the connections between the representation theory of the quantum Euclidean group and the so-called Hahn--Exton $q$-Bessel functions. In \cite{Noels}, Noels presented an algebraic construction of Fourier transforms between $\CO_q(E_2)$ and its Hopf dual. Brzezi\'{n}ski \cite{Brzezinski} constructed a complex of integral forms on the quantum Euclidean group based on twisted multi-derivations, and its isomorphism with the corresponding de Rham complex was established as well. Recently, in \cite{Aziziheris-Fakhri-Laheghi}, the authors construct $\CO_q(E_2)$ based on a new associative multiplication on $2\times 2$ matrices.  Recall that, as an abstract associative $\mK$-algebra, $\CO_q(E_2)$ is generated by $a, a^{-1}$, $b$, $c$ subject to the following relations
\begin{equation*}
ab=qba,  \quad ac=qca, \quad  bc=cb, \quad  aa^{-1}=a^{-1}a=1.
\end{equation*}
Clearly, $\CO_q(E_2)$ can be presented as a skew polynomial algebra $\CO_q(E_2)=\mK[b,c][a^{\pm 1}; \sigma]$ where $\sigma$ is the automorphism of the polynomial algebra $\mK[b,c]$ defined by $\s(b)=qb$ and $\s(c)=qc$. In particular, $\CO_q(E_2)$ is a Noetherian domain of Gelfand--Kirillov dimension 3.  The algebra $\CO_q(E_2)$ is made into a Hopf algebra by defining comultiplication, counit and antipode on the generators as follows
\begin{equation*}
\begin{aligned}
&\D(a)=a\otimes a ,   & \quad \varepsilon(a)&=1,  \quad &S(a)&=a^{-1},  \\
& \D(b)=b \otimes a^{-1}+a\otimes b, & \varepsilon(b)&=0, & S(b)&=-q^{-1}b, \\
& \D(c)=c\otimes a + a^{-1}\otimes c, & \varepsilon(c)&=0, & S(c)&=-qc. 
\end{aligned}
\end{equation*}

The Euclidean algebra $\mathfrak{e}_2$ is the complexification of the Lie algebra of the Euclidean group $E_2$. More precisely, $\mathfrak{e}_2$ is a 3-dimensional Lie algebra  with basis $\{ h, e,f \}$ and Lie bracket given by $[h,e]=2e$, $[h,f]=-2f$, and $[e,f]=0$. In particular, $\mathfrak{e}_2$ is a solvable Lie algebra. 
The quantized universal enveloping algebra $U_q(\mathfrak{e}_2)$ of the Lie algebra $\mathfrak{e}_2$ is the associative $\mK$-algebra generated by $K, K^{-1}$, $E$ and $F$ subject to the following relations
\begin{equation*}
KE=q^2EK, \;\;\;KF=q^{-2}FK, \;\;\;EF=FE, \;\;\;KK^{-1}=K^{-1}K=1. 
\end{equation*}
The defining relations of $U_q(\mathfrak{e}_2)$ can be obtained from those of $U_q(\sl_2)$ by taking certain limit, see e.g. \cite{Koelink} for this limit process. A different version of quantized universal enveloping algebra of $\mathfrak{e}_2$ based on the contraction procedure was given in \cite{Celeghini-Giachetti-Sorace-Tarlini}. 
Clearly, $U_q(\mathfrak{e}_2)$ can be presented as a skew polynomial algebra $U_q(\mathfrak{e}_2)=\mK[E,F][K^{\pm 1};\tau]$ where $\tau$ is the automorphism of the polynomial algebra $\mK[E,F]$ defined by $\tau(E)=q^2 E$ and $\tau(F)=q^{-2}F$. The algebra $U_q(\mathfrak{e}_2)$ is made into a Hopf algebra with structure maps given by
\begin{equation*}
\begin{aligned}
&\D(K)=K\otimes K ,  &\quad  \varepsilon(K)&=1,  \quad &S(K)&=K^{-1},  \\
& \D(E)=E\otimes K +1 \otimes E, & \varepsilon(E)&=0, & S(E)&=-EK^{-1}, \\
& \D(F)=F\otimes 1 + K^{-1}\otimes F, & \varepsilon(F)&=0, & S(F)&=-KF. 
\end{aligned}
\end{equation*}

There exists a non-degenerate dual pairing between $\CO_q(E_2)$ and $U_q(\mathfrak{e}_2)$. For the notion of a dual pairing of Hopf algebras, see for instance \cite[Subsection 1.2.5]{Klimyk-Schmudgen}.  More precisely, there is a unique non-degenerate dual pairing $\langle\cdot, \cdot \rangle : U_q(\mathfrak{e}_2) \times \CO_q(E_2) \rightarrow \mK$, with the non-vanishing pairs between the generators of  $U_q(\mathfrak{e}_2)$ and $\CO_q(E_2)$ as follows: $\langle K, a \rangle =q^{-1}$, $\langle K, a^{-1} \rangle =q$, $\langle E, c\rangle = \langle F, b \rangle =1$. This dual pairing turns the quantum Euclidean group $\CO_q(E_2)$ into a left $U_q(\mathfrak{e}_2)$-module algebra with the action determined by $u\cdot x =\sum x_{(1)}\langle u, x_{(2)}\rangle$ where $u\in U_q(\mathfrak{e}_2)$ and $x\in \CO_q(E_2)$. Here and below we use Sweedler's notation $\D(x)=\sum x_{(1)} \otimes x_{(2)}$ for the comultiplication of a Hopf algebra.  On generators, the action takes the following form
\begin{equation*}
\begin{aligned}
&K\cdot a =q^{-1}a, &K \cdot b&=qb, &K \cdot c&=q^{-1}c,\\
&E\cdot a =0, &E \cdot b&=0, &E \cdot c&=a^{-1},\\
&F\cdot a =0, &F \cdot b&=a, &F \cdot c&=0.
\end{aligned}
\end{equation*}
Therefore, one can form the smash product algebra $D_q(E_2):=\CO_q(E_2) \rtimes U_q(\mathfrak{e}_2)$, and we call $D_q(E_2)$ the \emph{Heisenberg double of the quantum Euclidean group} $\CO_q(E_2)$. The notion of smash product (also called cross product) has proved to be very useful in studying Hopf algebra actions, see \cite[Section 4.1]{Montgomery}. Heisenberg doubles are important smash product algebras, for the general construction of the Heisenberg double of a Hopf algebra, see \cite{Lu-HeisenbergDouble}. 
Note that $\CO_q(E_2)$ and $U_q(\mathfrak{e}_2)$ are naturally viewed as subalgebras of $D_q(E_2)$, and the cross commutation relations in $D_q(E_2)$ can be written as $ux=\sum (u_{(1)}\cdot x)u_{(2)}$ where $u\in U_q(\mathfrak{e}_2)$ and $x\in \CO_q(E_2)$. Inside the algebra $D_q(E_2)$, the left action of $U_q(\mathfrak{e}_2)$ on $\CO_q(E_2)$ can be nicely expressed by the left adjoint action of the Hopf algebra $U_q(\mathfrak{e}_2)$. One can think of $U_q(\mathfrak{e}_2)$ acts as generalized differential operators on the quantum Euclidean group. Therefore, the Heisenberg double $D_q(E_2)$ can be considered as an algebra of differential operators with coefficients in $\CO_q(E_2)$. As an abstract algebra, the defining relations of $D_q(E_2)$ are given below.
\begin{definition} \label{8Aug21}  
	The Heisenberg double $D_q(E_2)$ of the quantum Euclidean group  is the $\mK$-algebra generated by the elements $a, a^{-1}, b,c$, $K$, $K^{-1}$, $E$ and $F$ subject to the following defining relations
	\begin{equation*}
	\begin{aligned}
	ab&=qba,  &ac&=qca, &bc&=cb,  &aa^{-1}&=a^{-1}a=1, \\
	KE&=q^2EK, &KF&=q^{-2}FK, &EF&=FE, &KK^{-1}&=K^{-1}K=1, 
	\end{aligned}
	\end{equation*}
	\begin{equation*}
	\begin{aligned}
	&	Ka=q^{-1}aK, &\quad Kb&=qbK,& \quad Kc&=q^{-1}cK, && \\
	&Ea=aE, &Eb&=bE, &Ec&=cE+a^{-1}K, &&\\
	&	Fa=qaF, &Fb&=q^{-1}bF+a, &Fc&=qcF. &&
	\end{aligned}
	\end{equation*}
\end{definition}

To simplify notation, we abbreviate $\CO_q:=\CO_q(E_2)$, $U_q:=U_q(\mathfrak{e}_2)$, and $D_q:=D_q(E_2)$. 
The algebra $D_q$ will be the primary interest to us in this paper. We are concerned with ring-theoretic properties of $D_q$ and its representation theory. 

If $H$ is a finite-dimensional Hopf algebra and $H^*$ acts on $H$ via coregular action, then the Heisenberg double $D(H):=H \rtimes H^*$ is isomorphic to the algebra $\End_{\mK}(H)$, in particular, $D(H)$ is a simple algebra, see \cite[Corollary 9.4.3]{Montgomery}. The Heisenberg double of $\SL_2$, denoted $D(\SL_2)$, is the smash product of the coordinate ring $\CO(\SL_2)$ of $\SL_2$ with the universal enveloping algebra $U(\sl_2)$. Various classes of representations of $D(\SL_2)$ has been investigated in \cite{Tao-HeisenbergSL2}, and from there we know that $D(\SL_2)$ is a central simple algebra \cite[Proposition 3.1]{Tao-HeisenbergSL2}.  In contrast, the Heisenberg double $D_q$ of the quantum Euclidean group, studied in the present paper, has quite different properties. It is proved that the centre of $D_q$ is trivial (Proposition \ref{a15Oct21}), but $D_q$ is not a simple algebra. We give an explicit description of the prime, primitive and maximal ideals of $D_q$, together with the containments of prime ideals. All the prime factors of $D_q$ are explicitly described (Proposition \ref{a17Oct21}). It turns out that every prime factor of $D_q$ is a Noetherian domain that can be presented as a generalized Weyl algebra in the sense of \cite{Bavula-GGWA} (not in the usual sense of \cite{Bav-GWArep}). As a result, we obtain that $D_q$ has no finite-dimensional representations, and that $D_q$ cannot have a Hopf algebra structure (Corollary \ref{B19Oct21}).  The prime spectra of $\CO_q$ and $U_q$ are explicitly described as well (Theorem \ref{A16Oct21} and Theorem \ref{B16Oct21}). We see that all prime ideals of $\CO_q$ and $U_q$ are completely prime. 

There has been much interest in the study of automorphism groups of quantum algebras, see for instance \cite{Launois-Lenagan-PrimAutQuantMat,Launois-Lenagan-AutQuantMat,Yakimov-LauLenConj,Yakimov-AndDumConj} and references therein. 
It is a well-known fact that quantization leads to rigidity and so puts limits on the automorphism groups of quantum algebras. This means that the group of automorphisms of a quantum algebra is usually very ``small''. In Section \ref{AutGroup}, we determine the automorphism groups of the algebras $\CO_q$, $U_q$ and $D_q$. The method is to use the invariance of the set of height one prime ideals under the action of the automorphism group. Such method has been used successfully by Rigal \cite{Rigal} to calculate the automorphism group of the quantized Weyl algebra, and by Launois and Lenagan \cite{Launois-Lenagan-PrimAutQuantMat} to compute the automorphism group of the algebra $\CO_q(M_{m,n})$ of quantum matrices in the case where $m \neq n$.
We also describe the orbits of the action of the automorphism group on the prime spectrum.

A $D_q$-module $M$ is called a weight module if $M=\oplus_{\mu\in \mK^*}M_{\mu}$ where $M_{\mu}=\{x\in M\,|\,Kx=\mu x \}$. Note that each weight space $M_{\mu}$ is a module over the centralizer $\sC:=C(K,D_q)$ of the element $K$ in $D_q$.  The structure and representations of the algebra $\sC$ store a lot of information about weight $D_q$-modules. We find explicit generators and defining relations for the algebra $\sC$, and show that $D_q$ is free as a left and right $\sC$-module (Proposition \ref{a16Nov21}). The centre of $\sC$ is proved to be a Laurent polynomial algebra $\mK[K, K^{-1}]$.  For any $\chi \in \mK^*$, define $\sC(\chi):=\sC/\sC(K-\chi)$. We give a classification of all prime, primitive and maximal ideals of $\sC(\chi)$ (Theorem 5.5). It is interesting to note that all the prime factors of $\sC(\chi)$ can be presented explicitly as generalized Weyl algebras (Proposition \ref{d16Nov21}). From this we are able to give a classification of all simple unfaithful $\sC(\chi)$-modules. The next step is to investigate weight $D_q$-modules using the results obtained for the centralizer $\sC$. We show that for any  weight $D_q$-module, all its weight spaces are infinite-dimensional. For any algebra $R$, we denote by $\widehat{R}$ the set of isomorphism classes of simple $R$-modules. We show in Theorem \ref{B17Nov21} that the following map defines a surjection
\begin{equation*}
\widehat{\sC} \twoheadrightarrow \widehat{D_q}({\rm weight}), \quad M \mapsto \ind_{\sC}^{D_q}(M):=D_q \otimes_{\sC} M. 
\end{equation*}
In other words, any simple weight $D_q$-module is isomorphic to a module that induced from a simple $\sC$-module. Thus the problem of classification of simple weight $D_q$-modules is essentially reduced to a problem of classification of all simple modules over the algebra $\sC$. In particular, the classification of simple unfaithful $\sC(\chi)$-modules gives rise to a lot of simple weight $D_q$-modules. We also provide several examples of simple faithful $\sC(\chi)$-modules together with the description of the corresponding simple weight $D_q$-modules (Proposition \ref{20Nov21}). 

In the final section, we consider simple $D_q$-modules that are $\mK[a,a^{-1}]$-torsion. These are not weight modules as defined above. But notice that the element $a$ is normal in $D_q$ which plays a similar role to that of $K$, it is natural to call them \emph{$a$-weight modules}. We show that the problem of classification of simple $a$-weight modules is  essentially reduced to the problem of classification of all simple modules over the centralizer $C(a,D_q)$ of the element $a$ in $D_q$ (Theorem \ref{A22Nov21}). We give explicit generators and defining relations for $C(a,D_q)$, and show that the two algebras $C(a, D_q)$ and $C(K, D_q)$ are isomorphic. From this we can construct a lot of simple $a$-weight $D_q$-modules. 

The paper is organized as follows. In Section \ref{Preliminaries}, the centre of $D_q$ is determined, some centralizers are calculated which turn out to be quantum homogeneous spaces of $\CO_q$. In Section \ref{PrimeSpectra}, we describe the prime, completely prime, primitive and maximal spectra of the algebras $\CO_q$, $U_q$ and $D_q$. In Section \ref{AutGroup}, we determine the automorphism groups of the algebras $\CO_q$, $U_q$ and $D_q$. Section \ref{SimpleWeight} is devoted to the investigation of weight $D_q$-modules.  In Section \ref{aWeightMod}, we study $a$-weight $D_q$-modules. 

\section{Preliminaries} \label{Preliminaries} 

In this section, we prove that the centre of the Heisenberg double $D_q$ is trivial (Proposition \ref{a15Oct21}). 
Some centralizers are calculated (Lemma \ref{a13Nov21}), they provide examples of quantum homogeneous spaces for the quantum Euclidean group $\CO_q$.

We begin by recalling the notion of a quantum polynomial algebra. 
Let $Q=(q_{ij})$ be an $n \times n$ matrix whose entries $q_{ij}\in \mK^*$ satisfying $q_{ii}=q_{ij}q_{ji}=1$. Fix an integer $r$ such that $0 \leq r \leq n$, the \emph{quantum polynomial algebra}
$\Lambda:=\mK_Q[X_1^{\pm 1}, \ldots, X_r^{\pm 1}, X_{r+1}, \ldots, X_n]$
is the associated $\mK$-algebra generated by the indeterminates $X_1, \ldots, X_n$, $X_1^{-1},\ldots, X_r^{-1}$ subject to the following defining relations: $X_iX_j=q_{ij}X_j X_i$ ($1 \leq i, j \leq n$), and $X_i X_i ^{-1}=X_i^{-1}X_i=1$ ($i=1, \ldots, r$). Note that $\Lambda$ is the quantized coordinate ring of an affine space if $r=0$, and a quantum torus if $r=n$. 

For a ring $R$, we denote by $Z(R)$ the centre of $R$. An element $r\in R$ is \emph{normal} if $rR=Rr$.
Clearly, both $\CO_q$ and $U_q$ are examples of quantum polynomial algebras. It is easy to see that the centre of $\CO_q$ is trivial, and the centre of $U_q$ is a polynomial algebra in one variable $C:=EF$; that is  $Z(U_q)=\mK[C]$. The algebra $D_q$ can be expressed as an iterated Ore extension over the algebra $\CO_q$ of the  form
\begin{equation} \label{OreExt} 
D_q=\CO_q[K^{\pm 1};\s_1][E;\s_2, \d_2][F;\s_3,\d_3]. 
\end{equation}
In particular, $D_q$ is a Noetherian domain of Gelfand--Kirillov dimension 6. We note that the algebra $D_q$ has a PBW type basis, and the elements $K^{\pm 1}$ and $a^{\pm 1}$ are normal in $D_q$. 
Recall that an involution $*$ on an algebra $R$ is a $\mK$-algebra anti-automorphism such that $a^{**}=a$ for all $a \in R$. The algebra $D_q$ is involutive with an involution $*$ on $D_q$ defined by the following rule
\begin{equation} \label{Involution} 
\begin{aligned}
&a^*=a^{-1}, &b^*&=-qc, &c^*&=-q^{-1}b, \\
&K^*=K, &E^*&=qKF, &F^*&=qK^{-1}E. 
\end{aligned}
\end{equation}

Set
\begin{equation} \label{phipsi} 
\phi:=(1-q^2)Fb+q^2 a, \quad {\rm and} \quad  \psi:=(1-q^{-2})Ec+q^{-2}a^{-1}K.
\end{equation}
The elements $\phi$ and $\psi$ can be written more symmetrically as $\phi=Fb-qbF$ and $\psi=Ec-q^{-2}cE$. It is easily checked that $\psi=q^{-3}K \phi^*$ and $\phi=q^2K^{-1}\psi^*$. Moreover, $\phi$ and $\psi$ are normal elements of $D_q$, their commutation relations with the generators of $D_q$ are tabulated in the following
\begin{equation*}
\begin{aligned}
\phi K &=q K \phi, &\quad \phi a &= a\phi, &\quad \psi K &= q^{-1}K \psi, &\quad \psi a &= q^{-1}a \psi, \\
\phi E &=E\phi,  &\phi b&= q^{-1}b\phi,  &\psi E&= E \psi, &\psi b&= b\psi, \\
\phi F&=qF\phi,  &\phi c&=qc \phi,  &\psi F&=q^{-1}F \psi, &\psi c&=c\psi. 
\end{aligned}
\end{equation*}
Using the above relations, one easily verifies that $\phi \psi =q \psi \phi$. 
The following proposition shows that the centre of $D_q$ is trivial. 

\begin{proposition} \label{a15Oct21} 
	$Z(D_q)=\mK$. 	
\end{proposition}
\begin{proof}
	Let $D_qS^{-1}$ be the localization of $D_q$ at the Ore set $S=\{ b^i c^j \,|\, i,j\in \N \}$ of the algebra $D_q$. By the expressions of the elements $\phi$ and $\psi$ in (\ref{phipsi}), in the algebra $D_qS^{-1}$,  one can replace the generators $E$, $F$ by the elements $\psi$ and   $\phi$, respectively. It follows that the algebra $D_qS^{-1}$ can be presented as a quantum polynomial algebra of the following form 
	\begin{equation} \label{DqS}  
	D_qS^{-1}=\mK_Q[a^{\pm 1}, \,b^{\pm 1}, \,c^{\pm 1}, \,K^{\pm 1}, \,\psi, \,\phi]
	\end{equation}
	where $Q=(q_{ij})$ is a $6\times 6$ matrix whose entries are defined by $q_{ij}=q^{d_{ij}}$ wtih ${\bf D}=(d_{ij})$ given by 
	\begin{equation*}
	{\bf D}=\begin{psmallmatrix}
	\phantom{-}0 & \phantom{-}1& \phantom{-}1& \phantom{-}1& \phantom{-}1&\phantom{-}0\\
	-1& \phantom{-}0&\phantom{-}0&-1&\phantom{-}0&\phantom{-}1 \\
	-1 & \phantom{-}0 &\phantom{-}0 &\phantom{-}1 &\phantom{-}0&-1 \\
	-1 & \phantom{-}1&-1&\phantom{-}0 &\phantom{-}1 &-1\\
	-1& \phantom{-}0&\phantom{-}0&-1&\phantom{-}0&-1\\
	\phantom{-}0&-1&\phantom{-}1&\phantom{-}1&\phantom{-}1&\phantom{-}0
	\end{psmallmatrix}.
	\end{equation*}	
	An easy computation shows that $\ker({\bf D})=0$, and from this we conclude that $Z(D_qS^{-1})=\mK$. Therefore,  $Z(D_q)=D_q \cap Z(D_q S^{-1})=\mK$. 
\end{proof}

\begin{remark}
	Recall that a quantum Weyl field is a skew field of fraction of a quantum affine space. We say that a $\mK$-algebra $R$ admitting a skew field of fraction $\Frac(R)$ satisfies the \emph{quantum Gelfand--Kirillov conjecture} if $\Frac(R)$ is isomorphic to a quantum Weyl field over a purely transcendental field extension of $\mK$. Since $\Frac(D_q)=\Frac(D_qS^{-1})$, from (\ref{DqS}) it follows that the algebra $D_q$ satisfies the quantum Gelfand--Kirillov conjecture. 	
\end{remark}

For an element $r\in R$ and a subalgebra $S \subseteq R$, the subalgebra $C(r,S):=\{ s\in S\,|\, sr=rs \}$ is called the \emph{centralizer of $r$ with respect to $S$}. The following proposition describes the centralizers of the elements $K$, $E$ and $F$ with respect to the quantum Euclidean group $\CO_q \subset D_q$. It turns out that all these centralizers are quantum homogeneous spaces of $\CO_q$.

\begin{proposition} \label{a13Nov21} 
	Set $x:=ab$ and $y:=ca^{-1}$. 	
	\begin{enumerate}
		\item $C(K, \CO_q)=\Pi:=\mK\langle x, y\,|\, xy=q^2 yx \rangle$ is a quantum plane. Moreover, $\CO_q$ is free as a left and right $\Pi$-module with $\CO_q=\oplus_{i\in \Z} \Pi a^i =\oplus_{i\in \Z} a^i \Pi$.  
		\item $C(E, \CO_q)=\mK_q[a^{\pm 1}, b]$. 
		\item $C(F, \CO_q)=\mK[y]$ is a polynomial algebra. 
	\end{enumerate}
\end{proposition}
\begin{proof}
	1. Note that $\CO_q$ is a $\Z$-graded algebra with respect to the adjoint action of $K$, and the centralizer $C(K, \CO_q)$ consists of homogeneous elements of degree zero. More precisely, 
	\begin{equation*}
	C(K, \CO_q)=\bigoplus_{\stackrel{-i+j-k=0,}{i\in \Z, j,k\in \N}} \mK a^i b^j c^k=\bigoplus_{j,k\in \N}\mK a^{j-k}b^j c^k =\bigoplus_{j,k\in \N} \mK (ab)^j (ca^{-1})^k=\bigoplus_{j,k\in \N} \mK x^j y^k. 
	\end{equation*}
	Thus $C(K, \CO_q)$ is generated by the elements $x$ and $y$, and the elements $x^j y^k$ ($j,k\in \N$) form a basis of $C(K, \CO_q)$. Clearly, the relation $xy=q^2yx$ holds. Comparing bases yields that the algebras $C(K, \CO_q)$ and $\Pi$ are isomorphic.   The fact that $\CO_q$ is free as a left and right $\Pi$-module is evident. 
	
	2. Clearly $\mK_q[a^{\pm 1},b] \subseteq C(E,\CO_q)$, and $\CO_q=\oplus_{i\in \N} \mK_q[a^{\pm 1},b]c^i$.  Let $z=\sum_{i=0}^n p_i c^i$ where $p_i \in \mK_q[a^{\pm 1},b]$, be an element of $\CO_q$ that belongs to  $C(E, \CO_q)$. Using (\ref{Eci}) we obtain 
	\begin{equation*}
	Ez=zE+\sum_{i=0}^n \tfrac{1-q^{-2i}}{1-q^{-2}} p_i c^{i-1}a^{-1}K.
	\end{equation*}
	Since $Ez=zE$, we must have $n=0$, and hence $z\in \mK_q[a^{\pm 1},b]$. 
	
	3. It is evident that $\mK[y] \subseteq C(F, \CO_q)$, and that $\CO_q=\oplus_{i\in \N}\mK_q[a^{\pm 1},y]b^i$. Let $z=\sum_{i=0}^n p_i b^i$ where $p_i \in \mK_q[a^{\pm 1},y]$, be an element of $\CO_q$ that belongs to $C(F, \CO_q)$, and let $\s$ be the automorphism of $\mK_q[a^{\pm 1},y]$ defined by $\s(a)=qa$ and  $\s(y)=y$. Using (\ref{Fbi}) one deduces that 
	\begin{equation*}
	Fz=\sum_{i=0}^n \s(p_i)q^{-i}b^i F +\sum_{i=0}^n \tfrac{1-q^{-2i}}{1-q^{-2}}\s(p_i) ab^{i-1}. 
	\end{equation*}
	Since $Fz=zF$, comparing the coefficients of $F$ yields that $p_i = \s(p_i)q^{-i}$ and $\s(p_i)=0$ for $i \geq 1$. Hence $z$ belongs to the $\s$-invariants $\mK_q[a^{\pm 1},y]^{\s}$ in the algebra $\mK_q[a^{\pm 1},y]$. It is easy to see that $\mK_q[a^{\pm 1},y]^{\s}=\mK[y]$, and therefore, we have $C(F,\CO_q)=\mK[y]$ as desired.   
\end{proof}

The Laurent polynomial algebra $H=\mK[K,K^{-1}]$ is a Hopf algebra with $K$ a grouplike element. Recall that $\CO_q$ is a $\Z$-graded algebra with respect to the adjoint action of $K$.  For any homogeneous element $z\in \CO_q$, define $\rho: \CO_q \rightarrow \CO_q \otimes H$, $z \mapsto z \otimes K^{\deg z}$; the map $\rho$ makes the quantum Euclidean group $\CO_q$ into a right $H$-comodule algebra. It is evident that the algebra $\CO_q^{{\rm co} H}$ of coinvariants  coincides with the centralizer $\Pi=C(K, \CO_q)$. Moreover, using Proposition \ref{a13Nov21}.(1),  it is easy to show that the map $\beta: \CO_q \otimes_{\Pi} \CO_q \rightarrow \CO_q \otimes_{\mK}H$, $z\otimes z' \mapsto (z\otimes 1)\rho(z')$ is a bijection. It follows that the extension $\Pi=\CO_q^{{\rm co} H} \subset \CO_q$ is a right  $H$-Galois extension, see \cite[Definiton 8.1.1]{Montgomery}. 
Recall that the $q$-binomial coefficients are given by the formula
\begin{equation*}
\begin{bmatrix}
n\\m
\end{bmatrix}_q:=\frac{(q;q)_n}{(q;q)_m (q;q)_{n-m}} \quad {\rm where}\,\,(q;q)_i:=(1-q)(1-q^2) \cdots (1-q^i). 
\end{equation*}
A subspace $C$ of a Hopf algebra $H$ is called a \emph{left coideal} if $\D(C) \subset H \otimes C$. The following lemma shows that the quantum plane $\Pi$ is a quantum homogeneous space of $\CO_q$. For an introduction to quantum homogeneous space, we refer to \cite[Section 11.6]{Klimyk-Schmudgen}. 
\begin{lemma}
	The comultiplication $\D$ is given on a basis element $x^m y^n$ of $\Pi$ by
	\begin{equation*}
	\D(x^m y^n)=\sum_{i=0}^m \sum_{j=0}^n \begin{bmatrix}
	m\\i
	\end{bmatrix}_{q^2} \begin{bmatrix}
	n\\j
	\end{bmatrix}_{q^2} x^i a^{2m-2(i+j)}y^{n-j} \otimes x^{m-i}y^j. 
	\end{equation*}
	In particular, $\D(\Pi) \subset \CO_q \otimes \Pi$, and the quantum plane $\Pi$ is a left coideal subalgebra of $\CO_q$. 
\end{lemma}
\begin{proof}
	This follows from the $q$-binomial formula, see for instance \cite[Proposition 2.2]{Klimyk-Schmudgen}. 
\end{proof}

\begin{remark}
	\begin{enumerate}
		\item 	From Proposition \ref{a13Nov21}, it follows that	$C(E,\CO_q)$ is a Hopf subalgebra of $\CO_q$, and $C(F,\CO_q)$ is a left coideal subalgebra of $\CO_q$. Moreover, for each $x\in \{ K, E, F \}$, the algebra $\CO_q$ is free, and in particular faithfully flat, as a left (and right) $C(x, \CO_q)$-module.  Hence, each centralizer $C(x, \CO_q)$ is a quantum homogeneous space of $\CO_q$. 
		\item In the classical limit $q \rightarrow 1$, the Heisenberg double $D_q$  specializes to the second Weyl algebra $A_2$ over the Laurent polynomial ring $\mK[K^{\pm 1}, a^{\pm 1}]$. If we denote by $\mathscr{D}$ the limit of the algebra $D_q$ as $q\rightarrow 1$, then $Z(\mathscr{D})=\mK[K^{\pm 1}, a^{\pm 1}]$, and for any $\chi, \alpha \in \mK^*$, the factor algebra $\mathscr{D}/\mathscr{D}(K-\chi, a-\alpha)$ is isomorphic to the second Weyl algebra $A_2(\mK)$. 
	\end{enumerate}
\end{remark}

\section{The prime and primitive spectra of $D_q(E_2)$} \label{PrimeSpectra}  

The main aim of this section is to describe the prime spectrum of the algebra $D_q$ (Theorem \ref{A19Oct21}). All prime factors of $D_q$ are explicitly described (Proposition \ref{a17Oct21}). It turns out that every prime factor of $D_q$ is a Noetherian domain that can be presented as a generalized Weyl algebra in the sense of \cite{Bavula-GGWA}. As a consequence, we obtain the classification of completely prime, maximal and primitive ideals of $D_q$; moreover, we show that $D_q$ has no finite-dimensional representations, and that $D_q$ cannot have a Hopf algebra structure (Corollary \ref{B19Oct21}).  The prime spectra of $\CO_q$ and $U_q$ are explicitly described as well (Theorem \ref{A16Oct21} and Theorem \ref{B16Oct21}).

Recall that a proper ideal $\gp$ of a ring $R$ is called a \emph{prime ideal} whenever, for all ideals $\ga, \gb$ of $R$, if $\ga \gb \subseteq \gp$ then either $\ga \subseteq \gp$ or $\gb \subseteq \gp$. We denote by $\Spec(R)$ the set of all prime ideals of $R$, which is a partially ordered set with respective to inclusions of prime ideals.   An ideal $\gp$ of $R$ is said to be \emph{completely prime}  if $R/\gp$ is a domain. As is well known, a completely prime ideal is prime, but in general, the converse does not hold.  Given a subset $\CS$ of a ring $R$, we denote by $(\CS)$ the two-sided ideal of $R$ generated by $\CS$. 

The following theorem describes the prime spectrum of the algebra $\CO_q(E_2)$.
\begin{theorem} \label{A16Oct21}  
	$\Spec(\CO_q)=\{ (0), \,(b), \,(c), \,(b,\,c), \,(b,\,c, \,a-\alpha)\,|\, \alpha\in \mK^*  \}$. All possible containments of prime ideals are depicted in the following diagram
	
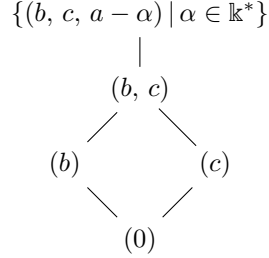
\begin{equation*}
	\begin{tikzpicture}[text centered]
	\node (bca) at (0,3)  {$\{(b,\,c,\,a-\alpha)\,|\,\alpha \in \mK^* \}$};
	\node (bc) at (0,2)  {$(b,\,c)$};
	\node (c) at (1,1)  {$(c)$};
	\node (b) at (-1,1)  {$(b)$};
	\node (0)  at (0,0)  {$(0)$};
	\draw [shorten <=-2pt, shorten >=-2pt] (bca) --(bc)-- (b)--(0)--(c)--(bc);
	\end{tikzpicture}
	\end{equation*}	
\end{theorem}

\begin{proof}
	Recall that $\CO_q=\mK[b,c][a^{\pm};\s]$ where $\s(b)=qb$ and $\s(c)=qc$. Clearly, the elements $b$ and $c$ are normal in $\CO_q$. The localization of $\CO_q$ at the Ore set $S=\{ b^i c^j \,|\, i,j\in \N \}$ is a quantum torus
	$\CO_qS^{-1}=\mK[b^{\pm 1}, c^{\pm 1}][a^{\pm 1};\s]$. An easy calculation shows that $Z(\CO_qS^{-1})=\mK$, then it follows from \cite[Corollary 1.5]{Goodearl-Letzter-QuantAffine}  that $\CO_qS^{-1}$ is  a simple algebra. Therefore, for any non-zero prime ideal $\gp$ of $\CO_q$, either $b \in \gp$ or $c\in \gp$. The prime ideals of $\CO_q$ that contain the element $b$ are in bijection with the prime ideals of the factor algebra $\CO_q/(b)\simeq \mK[c][a^{\pm 1};\s]$ where $\s(c)=qc$.  Notice that $\CO_q/(b)$ is isomorphic to a localization of the quantum plane, whose prime spectrum is easily obtained $\Spec(\CO_q/(b))=\{ (0), (c), (c, a-\alpha)\,|\, \alpha \in \mK^* \}$. Thus the set of prime ideals of $\CO_q$ that contain the element $b$ is $\{ (b), (b,c), (b,c,a-\alpha)\,|\, \alpha \in \mK^* \}$. Similarly,  the prime ideals of $\CO_q$ contain the element $c$ are in bijection with the prime ideals of the algebra $\CO_q/(c)\simeq \mK[b][a^{\pm 1};\s]$ where $\s(b)=qb$, which is isomorphic to $\CO_q/(b)$. It follows that the set of prime ideals of $\CO_q$ that contain the element $c$ is $\{ (c), (b,c), (b,c,a-\alpha)\,|\, \alpha \in \mK^* \}$. All the prime ideals of $\CO_q$ are shown in the diagram and the inclusions of prime ideals are obvious. This completes the proof. 
\end{proof}

The following theorem describes the prime spectrum of the algebra $U_q(\mathfrak{e}_2)$. 
\begin{theorem} \label{B16Oct21} 
	$\Spec(U_q)=\{ (0), (E), (F), (E,F), (E,F, K-\gamma)\,|\, \gamma \in \mK^*  \} \cup \{ (C-\beta)\,|\, \beta\in \mK^* \}$ where $C=EF$. All the containments of prime ideals of $U_q$ are shown in the following diagram
	
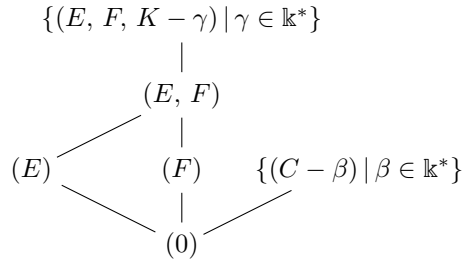
\begin{equation*}
	\begin{tikzpicture}[text centered]
	\node (EFK) at (0,3)  {$\{(E,\,F,\,K-\gamma)\,|\,\gamma \in \mK^* \}$};
	\node (EF) at (0,2)  {$(E,\,F)$};
	\node (F) at (0,1) {$(F)$};
	\node (C) at (2,1)  {\qquad$\{(C-\beta)\,|\, \beta \in \mK^*\}$};
	\node (E) at (-2,1)  {$(E)$};
	\node (0)  at (0,0)  {$(0)$};
	\draw [shorten <=-2pt, shorten >=-2pt] (EFK) --(EF)-- (F)--(0)--(E)--(EF);
	\draw [shorten <=-2pt, shorten >=-2pt] (C) --(0);
	\end{tikzpicture}
	\end{equation*}	
\end{theorem}
\begin{proof}
	Recall that the algebra $U_q=\mK[E,F][K^{\pm 1};\tau]$ where $\tau(E)=q^2E$ and $\tau(F)=q^{-2}F$, and that $Z(U_q)=\mK[C]$. Since the element $E$ is normal in $U_q$, the prime spectrum of $U_q$ decomposes into a disjoint union of the following two sets $S_1=\{ \gp \in \Spec(U_q)\,|\, E \in \gp \}$ and $S_2=\{ \gp \in \Spec(U_q)\,|\, E \notin \gp \}$. Clearly, the set $S_1$ is homeomorphic to the prime spectrum of the factor algebra $U_q/(E)\simeq \mK[F][K^{\pm 1};\tau]$ where $\tau(F)=q^{-2} F$. Since $U_q/(E)$ is a localization of the quantum plane, whose prime spectrum is easily described $\Spec(U_q/(E))=\{ (0), (F), (F, K-\gamma)\,|\, \gamma \in \mK^* \}$. It follows that $S_1=\{ (E), (E,F), (E,F, K-\gamma)\,|\, \gamma \in \mK^* \}$. The set $S_2$ is homeomorphic to the prime spectrum of the localization $U_q[E^{-1}]$ obtained by inverting the element $E$. Notice that $U_q[E^{-1}]$ can be presented as a tensor product of algebras
	\begin{equation*}
	U_q[E^{-1}] \simeq \mK[C] \otimes \mK[E^{\pm 1}][K^{\pm 1};\tau],
	\end{equation*}
	where the quantum torus $\mK[E^{\pm 1}][K^{\pm 1};\tau]$ is a simple algebra. It follows that the prime spectrum of $U_q[E^{-1}]$ is homeomorphic to the prime spectrum of the polynomial algebra $\mK[C]$. Consequently, $S_2=\{ (0)\}\cup \{\gm_{\beta}\cap U_q\,|\,\beta \in \mK \}$ where $\gm_{\beta}:=U_q[E^{-1}](C-\beta)$. For $\beta \in \mK^*$, it is easy to see that $U_q/(C-\beta) \simeq \mK[E^{\pm 1}][K^{\pm 1};\tau]$ is a simple algebra. In other words, the ideal $(C-\beta)$ of $U_q$ is maximal. Therefore in this case $\gm_{\beta}\cap U_q=(C-\beta)$. For $\beta=0$, note that $\gm_0=U_q[E^{-1}]C=U_q[E^{-1}]F$. Let us show that $\gm_0 \cap U_q=(F)$. The inclusion ``$\supseteq$'' is obvious. For the converse, if $x \in \gm_0 \cap U_q$ then $E^i x\in (F)$ for some positive integer $i$. Notice that $(F)$ is a completely prime ideal of $U_q$ and $E\notin (F)$, since $U_q/(F)\simeq \mK[E][K^{\pm 1};\tau]$ is a domain. We obtain that $x\in (F)$, and this proves $\gm_0 \cap U_q=(F)$. Therefore, $S_2=\{ (0) \}\cup \{ (F), (C-\beta)\,|\, \beta \in \mK^* \}$. All the prime ideals are shown in the diagram, and the inclusions are obvious. 	
\end{proof}

From the above description, we see that all prime ideals of $\CO_q$ and $U_q$ are completely prime. 

For an algebra $R$, we denote by $\widehat{R}$ the set of isomorphism classes of simple $R$-modules. If $\mathscr{P}$ is an isomorphism invariant property on simple $R$-modules, then $\widehat{R}(\mathscr{P})$ is the set of isomorphism classes of simple $R$-modules that satisfy property $\mathscr{P}$. 
The following corollary gives a classification of all simple finite-dimension $\CO_q$-modules and $U_q$-modules. 
\begin{corollary} \label{d16Oct21} 
	\begin{equation*}
	\widehat{\CO_q}\,({\rm fin. dim.})=\{ \CO_q/\gm_{\alpha} \,|\, \alpha \in \mK^*\},\,\,{\rm and}\,\,\,\, \widehat{U_q}\,({\rm fin. dim.})=\{ U_q/\gn_{\gamma} \,|\, \alpha \in \mK^*\}
	\end{equation*}
	where $\gm_{\alpha}=(b,c, a-\alpha)$ and $\gn_{\gamma}=(E,F, K-\gamma)$ are maximal ideals of $\CO_q$ and $U_q$, respectively. In particular, all these simple modules are of dimension one. 
\end{corollary}
\begin{proof}
	For a simple finite-dimensional $R$-module $M$,  its annihilator is a maximal ideal of $R$ with finite codimension, and $M$ is a simple module over the factor algebra $R/\ann_R(M)$. By Theorem \ref{A16Oct21}, the maximal ideals of $\CO_q$ are $\gm_{\alpha}$, $\alpha \in \mK^*$. By Theorem \ref{B16Oct21}, the maximal ideals of $U_q$ with finite codimension are $\gn_{\gamma}$, $\gamma \in \mK^*$. Notice  that $\CO_q/\gm_{\alpha} \simeq \mK$ and $U_q/\gn_{\gamma} \simeq \mK$, the result follows. 
\end{proof}

Let us now turn to the description of the prime spectrum of the algebra $D_q$. We will need the following simple observation.
\begin{lemma} \label{c16Oct21} 
	For all $i\in \N$, we have $(b^i)=D_q= (c^i)$. Moreover, $(b^i c^j)=D_q$ for all $i,j\in \N$. 
\end{lemma}
\begin{proof}
	In the algebra $D_q$, the following equalities hold for all $i \geq 1$:
	\begin{align}
	Fb^i &=q^{-i}b^i F+\big((1-q^{-2i})/(1-q^{-2})\big)ab^{i-1}, \label{Fbi} \\
	Ec^i&=c^i E +\big((1-q^{-2i})/(1-q^{-2})\big) c^{i-1}a^{-1}K. \label{Eci} 
	\end{align}
	The first equality can be easily verified by induction on $i$. The second equality follows from the first one by applying the involution $*$ defined in (\ref{Involution}). Since the elements $a$ and $K$ are invertible in $D_q$, the  equality (\ref{Fbi}) yields that $(b^i)=(b^{i-1})=\cdots =(b)=D_q$, and the equality (\ref{Eci}) implies that $(c^i)=(c^{i-1})=\cdots =(c)=D_q.$ Combining (\ref{Fbi}) with the identity $Fc=qcF$, we see that $b^{i-1}c^j \in (b^i c^j)$. Consequently, $(b^i c^j)=(c^j)=D_q$. 
\end{proof}

In the next proposition we shall give explicit descriptions of some factor algebras of $D_q$, which are in fact all the prime factors of $D_q$ (see Theorem \ref{A19Oct21}). It is interesting to note that all  these factor algebras belong to a large class of algebras which are  generalizations of the usual generalized Weyl algebras (GGWAs for short) introduced by Bavula in \cite{Bavula-GGWA}.  Let $R$ be a ring, $\s$ and $\tau$ be ring endomorphisms of $R$, and an element $\texttt{a}\in R$  be such that
\begin{equation} \label{GGWA1} 
\tau \sigma (\texttt{a})=\texttt{a}, \quad \texttt{a}r=\tau \sigma(r) \texttt{a}, \quad  \s(\texttt{a})r=\s \tau(r) \sigma(\texttt{a}), \quad {\rm for \,\,all}\,\, r\in R. 
\end{equation}
Recall from \cite{Bavula-GGWA} that, the GGWA $A:=R[x,y;\s, \tau, \texttt{a}]$ is a ring generated by $R$, $x$ and $y$ subject to the following defining relations:
\begin{equation*}
\begin{aligned}
xr&=\s(r)x, &yr&=\tau(r)y, \quad {\rm for \,\,all}\,\, r\in R, \\
yx&=\texttt{a}, &xy&=\s(\texttt{a}). 
\end{aligned}
\end{equation*}
The ring $R$ is called the base ring of the GGWA $A$. 
If $\s$ and $\tau$ are automorphisms and $\tau=\s^{-1}$, then the condition (\ref{GGWA1}) is equivalent to the fact that the element $\texttt{a}$ is central in $R$, and we obtain the usual definition of a generalized Weyl algebra (GWA for short), see for instance \cite{Bav-SimpleD[xy]Mod,Bav-GWArep}. In particular, every GWA is a GGWA. By \cite[Proposition 2.6]{Bavula-GGWA}, the GGWA is a domain if and only if the base ring $R$ is domain, $\texttt{a} \neq 0$ and $\s(\texttt{a}) \neq 0$.  

For a set $I=\{ i_1<\cdots< i_t \} \subseteq \{1,2,\ldots, 6\}$, we denote by $Q_I$ the $(6-t)\times (6-t)$ matrix obtained from $Q$ (defined in (\ref{DqS})) by deleting the rows $i_1, \ldots, i_t$ and columns $i_1, \ldots, i_t$. For example, $Q_{36}$ is the $4\times 4$ matrix obtained from $Q$ by deleting the third and sixth rows and columns.  The following proposition describes the prime factors of $D_q$, it also provides examples of GGWAs which are not the usual GWAs. 
\begin{proposition} \label{a17Oct21} 
	\
	\begin{enumerate}
		\item The factor algebra $D_q/D_q\phi$ is a domain that can be presented as a GGWA of the form
		\begin{equation*}
		D_q/D_q \phi \simeq \mK_{Q_{36}}[a^{\pm 1},\, b^{\pm 1}, \,K^{\pm 1}, \,\psi] \big[E,c;\s, \tau, \emph{\texttt{a}}\big], \quad \emph{\texttt{a}}=(1-q^{-2})^{-1}(\psi-a^{-1}K),
		\end{equation*} 
		where $\s$ and $\tau$ are automorphisms of the base ring defined by the following rule: $\s(a)=a$, $\s(b)=b$, $\s(K)=q^{-2}K$, $\s(\psi)=\psi$; and $	\tau(a)=q^{-1}a$, $\tau(b)=b$, $\tau(K)=qK$, $\tau(\psi)=\psi.$ Moreover, the centre of $D_q/D_q\phi$ is a polynomial algebra $\mK[\Theta]$ where $\Theta:=\psi b^{-1}$.

		\item The factor algebra $D_q/D_q\psi$ is a domain that can be presented as a GGWA of the form
		\begin{equation*}
		D_q/D_q \psi \simeq \mK_{Q_{25}}[a^{\pm 1}, c^{\pm 1}, K^{\pm 1}, \phi] \big[F,b;\s, \tau, \emph{\texttt{a}}\big], \quad \emph{\texttt{a}}=(q^{-1}-q)^{-1}(\phi-a),
		\end{equation*} 
		where $\s$ and $\tau$ are the automorphisms of the base ring defined by the rule: $\s(a)=qa$, $\s(c)=qc$, $\s(K)=q^2K$, $\s(\phi)=q^{-1}\phi$; and $\tau(a)=q^{-1}a$, $\tau(c)=c$, $\tau(K)=q^{-1}K$, $\tau(\phi)=q\phi$. Moreover, the centre of $D_q/D_q\psi$ is a polynomial algebra $\mK[\Omega]$ where $\Omega:=\phi K c^{-1}$.

		\item For any $\alpha \in \mK$, the factor algebra $D_q/(\phi, \,\psi-\alpha b)$ is a simple domain, which can be presented as a GGWA of the following form
		\begin{equation*}
		D_q/(\phi, \psi-\alpha b) \simeq \mK_{Q_{356}}[a^{\pm 1}, b^{\pm 1}, K^{\pm 1}] \big[E,c;\s, \tau,  \emph{\texttt{a}}\big], \quad \emph{\texttt{a}}=(1-q^{-2})^{-1}(\alpha b-a^{-1}K),
		\end{equation*}
		where $\s$ and $\tau$ are automorphisms of the base ring defined by the following rule: $\s(a)=a$, $\s(b)=b$, $\s(K)=q^{-2}K$; and $\tau(a)=q^{-1}a$, $\tau(b)=b$, $\tau(K)=qK.$ In particular, if $\alpha=0$ then $D_q/D_q(\phi,\psi)$ is a simple quantum torus
		\begin{equation*}
		D_q/D_q(\phi,\psi) \simeq \mathbb{T}:=\mK_{Q_{56}}[a^{\pm 1}, b^{\pm 1}, c^{\pm 1}, K^{\pm 1}]. 
		\end{equation*}
		
		\item For any $\beta \in \mK$, the factor algebra $D_q/(\phi-\beta c K^{-1}, \, \psi)$ is a simple domain, which can be presented as a GGWA of the following form
		\begin{equation*}
		D_q/(\phi-\beta c K^{-1}, \, \psi) \simeq \mK_{Q_{256}}[a^{\pm 1}, c^{\pm 1}, K^{\pm 1}] \big[F,b;\s, \tau, \emph{\texttt{a}}\big], \,\, \emph{\texttt{a}}=(q^{-1}-q)^{-1}(\beta cK^{-1}-a),
		\end{equation*} 
		where $\s$ and $\tau$ are the automorphisms of the base ring defined by the rule: $\s(a)=qa$, $\s(c)=qc$, $\s(K)=q^2K$; and $\tau(a)=q^{-1}a$, $\tau(c)=c$, $\tau(K)=q^{-1}K$. 
	\end{enumerate}
	
\end{proposition}
\begin{proof}
	1. Recall form (\ref{OreExt}) that the algebra $D_q$ can be presented as an iterated Ore extension of the form $D_q=A[F;\s_3, \d_3]$ where $A=\CO_q[K^{\pm1};\s_1][E;\s_2, \d_2]$. From the expression of $\phi$ and the fact that $a$ is invertible in $D_q$, one sees that the element $b$ is invertible in $D_q/D_q\phi$. It follows that $D_q/D_q\phi \simeq D_q[b^{-1}]/D_q[b^{-1}]\phi$ where $D_q[b^{-1}]$ is the localization of $D_q$ at the powers of the element $b$. From the expression of $\phi$, it is easily seen that $D_q[b^{-1}]/D_q[b^{-1}]\phi$ is isomorphic to the localization $A[b^{-1}]$ of the algebra $A$ at the powers of $b$. Clearly, the algebra $A[b^{-1}]$ can be presented as a GGWA as in the statement. In particular, $D_q/D_q\phi$ is a domain. It is straightforward to verify that $\Theta=\psi b^{-1}$ is a central element of $D_q/D_q\phi$. Considering the localization of the GGWA $D_q/D_q\phi$ at the powers of the element $c$, we obtain
	\begin{equation} \label{Dqphic} 
	(D_q/D_q\phi)[c^{-1}]\simeq \mK[\Theta] \otimes \mT
	\end{equation} 
	where the quantum torus $\mT=\mK_{Q_{56}}[a^{\pm 1}, b^{\pm 1}, c^{\pm 1}, K^{\pm 1}]$ is a central simple algebra. It follows that the centre of $(D_q/D_q\phi)[c^{-1}]$ is the polynomial algebra $\mK[\Theta]$, and therefore, $Z(D_q/D_q\phi)=\mK[\Theta]$. 
	
	2. By applying the involution $*$, see (\ref{Involution}), we see that the ideal $(D_q \phi)^*=D_q \psi$. Therefore, we have $D_q/D_q\psi \simeq (D_q/D_q\phi)^*$. Statement 2 then follows from statement 1. Notice that the element $\Theta^*=-q^{-4}\Omega$, so the centre of $D_q/D_q \psi$  is the polynomial algebra $\mK[\Omega]$. 
	
	3. The expression of $\CR:=D_q/(\phi, \psi- \alpha b)$ as a GGWA follows directly from statement 1. In particular, $\CR$ is a domain. If $\alpha=0$ then $\texttt{a}=-(1-q^{-2})^{-1}a^{-1}K$ is an invertible element, and therefore $c$ and $E$ are invertible in $D_q/(\phi, \psi)$. It follows that $D_q/(\phi,\psi) \simeq \mT$ is a quantum torus, which is a central simple algebra. For $\alpha \in \mK^*$, considering the localization $\CR[c^{-1}]$ of $\CR$ at the powers of the element $c$, one obtains $\CR[c^{-1}] \simeq \mT$. In particular, $\CR[c^{-1}]$ is a central simple algebra. Hence the algebra $\CR$ is simple if and only if $(\phi, \psi-\alpha b, c^i)=D_q$ for all $i \geq 1$, but this follows from Lemma \ref{c16Oct21} immediately. 
	
	4. The proof of statement 4 is similar to that of statement 3. 
\end{proof}

We are now in the position to give a classification of prime ideals of $D_q$. The generators of each prime ideal and inclusions between prime ideals are described explicitly.
\begin{theorem} \label{A19Oct21} 
	The prime spectrum of the algebra $D_q$ is given below,
	\begin{equation*}
	\Spec(D_q)=\{ (0), (\phi), (\psi), (\phi, \psi)  \} \,\sqcup \{  (\phi, \psi-\alpha b)\,|\, \alpha \in \mK^* \} \,\sqcup \,\{ (\phi-\beta c K^{-1}, \psi)\,|\, \beta\in \mK^* \}.
	\end{equation*}
	All the containments of the prime ideals of $D_q$ are shown in the following diagram
	\begin{equation*}
	\begin{tikzpicture}[text centered]
	\node (phisib) at (2,2) {\qquad\qquad $\{(\phi-\beta c K^{-1},\psi)|\beta \in \mK^*\}$};
	\node (phisia) at (-2,2) {$\{(\phi,\psi-\alpha b)|\alpha \in \mK^*\}$\qquad\qquad\qquad};
	\node (phisi) at (0,2) {$(\phi,\psi)$};
	\node (psi) at (1,1)  {$(\psi)$};
	\node (phi)  at (-1,1)  {$(\phi)$};
	\node (0)  at (0,0)  {$(0)$};
	\draw [shorten <=-2pt, shorten >=-2pt] (phisia) -- (phi) -- (0) -- (psi)--(phisib);
	\draw [shorten <=-2pt, shorten >=-2pt] (phi)--(phisi)--(psi);
	\end{tikzpicture} \label{SpL0} 
	\end{equation*} 	
\end{theorem}
\begin{proof}
	From Lemma \ref{c16Oct21}, it follows that any prime ideal of $D_q$ has no intersection with the Ore set $S=\{ b^i c^j \,|\, i,j\in\N \}$. Therefore, the prime spectrum of $D_q$ is homeomorphic to the prime spectrum of the algebra $D_qS^{-1}$. Now we turn to the description of $\Spec(D_qS^{-1})$.     If we denote by $D_q \mS^{-1}$ the localization of $D_q$ at the Ore set $\mS=\{ \mK^* b^i c^j \phi^k \psi^{\ell}\,|\,i,j,k,\ell \in \N \}$, then $D_q \mS^{-1}$ is actually the localization of $D_qS^{-1}$ obtained by inverting the elements $\phi$ and $\psi$. From (\ref{DqS}) it follows that 
	\begin{equation*}
	D_q \mS^{-1}=\mK_Q[a^{\pm 1}, b^{\pm 1}, c^{\pm 1}, K^{\pm 1}, \psi^{\pm 1}, \phi^{\pm 1}]
	\end{equation*}
	is a quantum torus where the matrix $Q$ is defined as in (\ref{DqS}). Notice that $Z(D_q\mS^{-1})=\mK$, we conclude that $D_q \mS^{-1}$ is a central simple algebra. Consequently, for any non-zero prime ideal $\gp$ of the algebra $D_qS^{-1}$, either $\phi \in \gp$ or $\psi \in \gp$. Thus the set of non-zero prime ideals of $D_qS^{-1}$ is decomposed into the union of the following two subsets: $\mathscr{S}(\phi):=\{ \gp \in \Spec(D_qS^{-1})\,|\, \phi \in \gp \}$ and $\mathscr{S}(\psi):=\{ \gp \in \Spec(D_qS^{-1})\,|\, \psi \in \gp \}$. Using the involution (\ref{Involution}), we see that $\mathscr{S}(\psi)=\mathscr{S}(\phi)^*$. Thus it suffices to describe the set $\mathscr{S}(\phi)$ of prime ideals of $D_qS^{-1}$ that contain the element $\phi$.
	Clearly, the set $\mathscr{S}(\phi)$ is homeomorphic to the prime spectrum of the factor algebra $D_qS^{-1}/D_qS^{-1}\phi$. From (\ref{DqS}), we obtain that  $D_qS^{-1}/D_qS^{-1}\phi  \simeq \mK[\Theta] \otimes \mT$ where $\Theta=\psi b^{-1}$, see also (\ref{Dqphic}). Thus the prime spectrum of $D_qS^{-1}/D_qS^{-1}\phi$ is homeomorphic to the prime spectrum of the polynomial algebra $\mK[\Theta]$, since $\mT$ is a central simple algebra. For a subset $\CG$ of $D_qS^{-1}$, let us denote by $(\CG)_S$ the two-sided ideal of $D_qS^{-1}$ generated by $\CG$. Then $\mathscr{S}(\phi)=\{ (\phi)_S, (\phi, \Theta-\alpha)_S\,|\, \alpha \in \mK \}$, since the field $\mK$ is algebraically closed. Applying the involution (\ref{Involution}), we obtain $\mathscr{S}(\psi)=\{ (\psi)_S, (\psi, \Omega-\beta)_S\,|\, \beta \in \mK \}$. We thus obtain the prime spectrum of the algebra $D_qS^{-1}$. 
	Now $\Spec(D_q)=\{ \gp \cap D_q\,|\, \gp \in \Spec(D_qS^{-1}) \}$, and we have to describe the generators of the prime ideal $\gp \cap D_q$ of $D_q$ for each $\gp \in \Spec(D_qS^{-1})$. 
	
	Clearly, $(0)_S \cap D_q=(0)$. Let us now show that $(\phi)_S \cap D_q=(\phi)$. The inclusion ``$\supseteq$'' is obvious. For the converse, if $x\in (\phi)_S\cap D_q$ then $xb^i c^j \in (\phi)$ for some $i,j\in \N$. From Proposition \ref{a17Oct21}.(1), the ideal $(\phi)$ of $D_q$ is completely prime and $b, c \notin (\phi)$. Thus $x \in (\phi)$ and we have $(\phi)_S \cap D_q=(\phi)$. Similarly, using Proposition \ref{a17Oct21}.(2), one can show that $(\psi)_S \cap D_q=(\psi)$. Now we prove that $(\phi, \Theta-\alpha)_S \cap D_q=(\phi, \psi-\alpha b)$. Notice that $(\phi, \Theta-\alpha)_S \cap D_q=(\phi, \psi-\alpha b)_S \cap D_q$ which contains the ideal $(\phi, \psi-\alpha b)$ of $D_q$. It must be an equality since $(\phi, \psi-\alpha b)$ is a maximal ideal of $D_q$, see Proposition \ref{a17Oct21}.(3). Similarly, using Proposition \ref{a17Oct21}.(4),  we obtain $(\psi, \Omega-\beta)_S \cap D_q=(\psi, \phi-\beta cK^{-1})$. This completes the description of $\Spec(D_q)$. All prime ideals of $D_q$ are shown in the diagram, and the containments are obvious. 
\end{proof}

From the diagram in Theorem \ref{A19Oct21}, we see that the prime spectrum of $D_q$ is homeomorphic to the prime spectrum of the quantum plane. Recall that, an ideal $\gp$ of a ring $R$ is called a \emph{primitive ideal} if $\gp$ is the annihilator of some simple $R$-module. The set $\Prim(R)$ of all primitive ideals is called the primitive spectrum of $R$. It is well-known that all primitive ideals of $R$ are prime, and all maximal ideals of $R$ are primitive. A prime ideal $\gp$ is said to be \emph{locally closed} in $\Spec(R)$ if the intersection of the prime ideals that properly contain $\gp$ is strictly larger than $\gp$. 
The following corollary describes the completely prime, maximal, and primitive ideals of the algebra $D_q$. As a consequence, $D_q$ has no finite-dimensional representations, and furthermore,  $D_q$ cannot have a Hopf algebra structure. 
\begin{corollary} \label{B19Oct21} 
	\
	\begin{enumerate}
		\item All prime ideals of $D_q$ are completely prime. 
		\item $\Max(D_q)=\{ (\phi, \psi-\alpha b)\,|\, \alpha \in \mK \} \,\sqcup\, \{ (\phi-\beta c K^{-1}, \psi)\,|\, \beta \in \mK^* \}$. 
		\item $\Prim(D_q)=\{(0) \} \,\sqcup\, \Max(D_q)$. 
		\item The Heisenberg double $D_q$ has no finite-dimensional modules. 
		\item The Heisenberg double $D_q$ cannot have a Hopf algebra structure. 
	\end{enumerate}
\end{corollary}
\begin{proof}
	1. All prime factor algebras of $D_q$ are described in Proposition \ref{a17Oct21}, which turn out to be domains. It means that all prime ideals of $D_q$ are completely prime.
	
	2. Statement 2 follows immediately from Theorem \ref{A19Oct21}. 
	
	3. First, maximal ideals of $D_q$ are primitive. Second, the ideals $(\phi)$ and $(\psi)$ of $D_q$ are not primitive, as the factor algebras $D_q/(\phi)$ and $D_q/(\psi)$ have non-trivial centres, see Proposition \ref{a17Oct21}.(1) and (2).  Since $D_q$ is a Jacobson ring, every locally closed prime ideal of $D_q$ is primitive. From the diagram in Theorem \ref{A19Oct21} we see that the ideal $(0)$ is locally closed, and thus primitive. 
	
	4. If $D_q$ has a finite-dimensional module, then $D_q$ has a finite-dimensional \emph{simple} module $M$. The annihilator of $M$ is then a primitive ideal of $D_q$ with finite codimension.  However, from Proposition \ref{a17Oct21}, all the primitive factors of $D_q$ are of infinite dimension, a contradiction. 
	
	5. If $D_q$ has a Hopf algebra structure, 	then the kernel of the counit $\varepsilon: D_q \rightarrow \mK$ is a maximal ideal of $D_q$ with codimension one. However, from statement 2 and Proposition \ref{a17Oct21}.(3) and (4), all maximal ideals of $D_q$ are of infinite codimension, a contradiction. 
\end{proof}

\section{The automorphism groups} \label{AutGroup} 

In this section, we investigate the groups of automorphisms of the algebras $\CO_q$, $U_q$ and $D_q$. The method is to use the invariance of the set of height one primes under the action of the automorphism group.  Using the classification of prime ideals of these algebras obtained in the previous section, explicit descriptions of their automorphism groups are obtained. The orbits of the action of the automorphism group on the prime spectrum are described. 

The following theorem describes the automorphism group of the algebra $\CO_q$. 
\begin{theorem} \label{23Oct21} 
	\
	\begin{enumerate}
		\item Let $\tau \in \Aut\, \CO_q$ be defined by $ \tau:$\,\, $a \mapsto a$, $b \mapsto c$, $c \mapsto b.$
		Then $\{ {\rm id}, \tau \}$ is a subgroup of $\Aut\,\CO_q$ isomorphic to the cyclic group $\mathbb{Z}_2$. 
		
		\item For any $i\in \Z$, there exists $\xi_i \in \Aut\,\CO_q$ such that $\xi_i:$\,\,  $a \mapsto a$,   $b \mapsto a^ib$, $c \mapsto a^i c.$
		Moreover, $\{ \xi_i \,|\, i \in \Z \}$ is a subgroup of $\Aut\,\CO_q$ isomorphic to $\Z$. 
		\item For all $(\alpha, \beta, \gamma) \in (\mK^*)^3$, there exists $\eta_{\alpha, \beta, \gamma} \in \Aut\,\CO_q$ such that 
		$\eta_{\alpha, \beta, \gamma}:$\,\, $a \mapsto \alpha a$, $b \mapsto \beta b$, $c \mapsto \gamma c$. 
		Moreover, $\{ \eta_{\alpha, \beta, \gamma}\,|\, (\alpha, \beta, \gamma)\in (\mK^*)^3 \}$ is a subgroup of $\Aut\,\CO_q$ isomorphic to the torus $(\mK^*)^3$. 
		
		\item $\Aut\,\CO_q  \simeq (\mK^*)^3 \rtimes (\Z \times \Z_2)$. 
	\end{enumerate}
\end{theorem}
\begin{proof}
	Statements 1--3 are clear. The subgroup $G$ of $\Aut\,\CO_q$ generated by the automorphisms $\tau$, $\xi_i$ ($i\in \Z$) and $\eta_{\alpha, \beta, \gamma}$ ($\alpha, \beta,\gamma \in \mK^*$)  is isomorphic to the semi-direct product $(\mK^*)^3 \rtimes (\Z \times \Z_2)$. It remains to show that $G=\Aut\,\CO_q$.  Let $\rho \in \Aut\,\CO_q$, we have to show that $\rho \in G$. From the diagram in Theorem \ref{A16Oct21}, the height one primes of $\CO_q$ are $(b)$ and $(c)$. So there are two options, either  the ideals $(b)$ and $(c)$ are $\rho$-invariant, or otherwise, they are interchanged. 
	
	If $(b)$ and $(c)$ are $\rho$-invariant, then $\rho(b)=xb$ and $\rho(c)=yc$ for some units $x$, $y$ of the algebra $\CO_q$. Note that the group $U(\CO_q)$ of units of the algebra $\CO_q$ is equal to $\{\mK^* a^i \,|\,i\in\Z  \}$. Therefore, $\rho(b)=\beta a^i b$ and $\rho(c)=\gamma a^j c$ for some $\beta, \gamma \in \mK^*$, $i,j\in \Z$. Using the identity $\rho(b)\rho(c)=\rho(c) \rho(b)$, we deduce that $i=j$. The automorphism $\rho$ preserves the group $U(\CO_q)$ of units. Thus either $\rho(a)=\alpha a$,  or otherwise, $\rho (a)=\alpha a^{-1}$ for some $\alpha \in \mK^*$. Let us show that the second case is impossible. Otherwise the equality $\rho(a)\rho(b)=q\rho(b) \rho(a)$ yields that $\alpha a^{-1} \cdot \beta a^i b= q\beta a^i b \cdot \alpha a^{-1}=q^2 \beta \alpha a^{i-1}b$. Hence $q^2=1$, a contradiction. It is easily verified that the map $\rho:$\,\,$ a \mapsto \alpha a,$  $b \mapsto \beta a^i b$,  $c \mapsto \gamma a^ic$
	defines an automorphism of $\CO_q$, and that $\rho=\xi_i \circ \eta_{\alpha, \beta, \gamma}$ belongs to $G$. 
	
	If the ideals $(b)$ and $(c)$ are interchanged under the action of $\rho$, then $\rho(b)=\beta a^i c$ and $\rho(c)=\gamma a^j b$ for some $\beta, \gamma \in \mK^*$, $i,j\in \Z$. The equality $\rho(b)\rho(c)=\rho(c)\rho(b)$ yields that $i=j$. The automorphism $\rho$ preserves the group $U(\CO_q)$ of units. Thus either $\rho(a)=\alpha a$,  or otherwise, $\rho (a)=\alpha a^{-1}$ for some $\alpha \in \mK^*$. The second case is impossible, since otherwise the equality $\rho(a)\rho(b)=q\rho(b) \rho(a)$ would imply that $q^2=1$. It is easily verified that the map $\rho:$\,\,$ a \mapsto \alpha a,$  $b \mapsto \beta a^i c$,  $c \mapsto \gamma a^ib$ defines an automorphism of $\CO_q$, and that $\rho =  \xi_i \circ \tau \circ \eta_{\alpha,\beta, \gamma}$ belongs to $G$. This completes the proof. 	
\end{proof}

\begin{corollary}
	The action of $\Aut\,\CO_q$ on $\Spec(\CO_q)$ has exactly 4 orbits:
	\begin{equation*}
	\{(0)\}, \quad \{ (b), \,(c) \},\quad \{ (b,c) \}, \quad \{ (b,c,a-\alpha)\,|\, \alpha\in \mK^* \}. 
	\end{equation*}
\end{corollary}
\begin{proof}
	The result follows from Theorem \ref{A16Oct21} and Theorem \ref{23Oct21} immediately. 
\end{proof}

The following theorem describes the group of automorphisms of the algebra $U_q$. 
It turns out that the automorphism groups of $U_q$ and $\CO_q$ are isomorphic. The proof is similar to that of Theorem \ref{23Oct21}; however, we include it for the convenience of the reader. 
\begin{theorem} \label{A24Oct21} 
	\
	\begin{enumerate}
		\item Let $\sigma \in \Aut\, U_q$ be defined by $\sigma:$\,\, $K \mapsto K^{-1}$, $E \mapsto F$, $F \mapsto E.$
		Then	$\{ {\rm id}, \sigma \}$ is a subgroup of $\Aut\,U_q$ isomorphic to the cyclic group $\mathbb{Z}_2$. 
		
		\item For any $i\in \Z$, there exists $\xi_i \in \Aut\,\CO_q$ such that
		$\xi_i:$\,\, $K \mapsto K$,  $E \mapsto K^i E$, $F \mapsto K^{-i} F$.
		Moreover, $\{ \xi_i \,|\, i \in \Z \}$ is a subgroup of $\Aut\,U_q$ isomorphic to $\Z$. 
		
		\item For all $(\alpha, \beta, \gamma) \in (\mK^*)^3$, there exists $\eta_{\alpha, \beta, \gamma} \in \Aut\,U_q$ such that $\eta_{\alpha, \beta, \gamma}:$\,\, $K \mapsto \alpha K$, $E \mapsto \beta E$, $F \mapsto \gamma F$.
		Moreover,  $\{ \eta_{\alpha, \beta, \gamma}\,|\, (\alpha, \beta, \gamma)\in (\mK^*)^3 \}$ is a subgroup of $\Aut\,U_q$ isomorphic to the torus $(\mK^*)^3$. 
		
		\item $\Aut\,U_q  \simeq (\mK^*)^3 \rtimes (\Z \times \Z_2)$. 
	\end{enumerate}
\end{theorem}
\begin{proof}
	Statements 1--3 are clear. The subgroup $G$ of $\Aut\,U_q$ generated by the automorphisms $\sigma$, $\xi_i$ ($i\in \Z$) and $\eta_{\alpha, \beta, \gamma}$ ($\alpha, \beta,\gamma \in \mK^*$)  is isomorphic to the semi-direct product $(\mK^*)^3 \rtimes (\Z \times \Z_2)$. It remains to show that $G=\Aut\,U_q$.  Let $\rho \in \Aut\,U_q$, we have to show that $\rho \in G$. 
	Recall that $Z(U_q)=\mK[C]$ where $C=EF$. Since automorphisms of $U_q$ preserve its centre, 	
	from the diagram in Theorem \ref{B16Oct21}, for any automorphism $\rho$ there are two options, either  the ideals $(E)$ and $(F)$ are $\rho$-invariant, or otherwise, they are interchanged. Note that the group of units of the algebra $U_q$ is equal to $\{\mK^* K^i \,|\,i\in\Z  \}$. 
	
	If $(E)$ and $(F)$ are $\rho$-invariant, then  $\rho(E)=\beta K^i E$ and $\rho(F)=\gamma K^j F$ for some $\beta, \gamma \in \mK^*$, $i,j\in \Z$. Using the identity $\rho(E)\rho(F)=\rho(F) \rho(E)$, we deduce that $i+j=0$. Since $\rho$ preserves the group of units, we have either $\rho(K)=\alpha K$,  or $\rho (K)=\alpha K^{-1}$ for some $\alpha \in \mK^*$. The second case is impossible, since otherwise the identity $\rho(K)\rho(E)=q^2 \rho(E) \rho(K)$ would imply that $q^4=1$. Now it is easily checked that the map $\rho:$\,\,$K \mapsto \alpha K$, $E \mapsto \beta K^i E$, $F \mapsto \gamma K^{-i}F$ defines an automorphism of $U_q$, and that $\rho=\xi_i \circ \eta_{\alpha,\beta,\gamma}$ belongs to $G$. 
	
	If $(E)$ and $(F)$ are interchanged under the action of $\rho$, then  $\rho(E)=\beta K^i F$ and $\rho(F)=\gamma K^j E$ for some $\beta, \gamma \in \mK^*$, $i,j\in \Z$. Then the equality $\rho(E) \rho(F)=\rho(F) \rho(E)$ yields that $i+j=0$. Since $\rho$ preserves the units of $U_q$, we have either $\rho(K)=\alpha K$,  or $\rho (K)=\alpha K^{-1}$ for some $\alpha \in \mK^*$. The first case is impossible, since otherwise the identity $\rho(K) \rho(E)=q^2 \rho(E)\rho(K)$ would imply that $q^4=1$. It is easy to verify that the map $\rho:$\,\,$K \mapsto \alpha K^{-1}$, $E \mapsto \beta K^i F$, $F \mapsto \gamma K^{-i}E$ defines an automorphism of $U_q$, and that $\rho=\xi_{-i} \circ \sigma \circ \eta_{\alpha,\beta,\gamma}$ belongs to $G$. This completes the proof. 
\end{proof}

\begin{corollary}
	The action of $\Aut\,U_q$ on $\Spec(U_q)$ has exactly 5 orbits:
	\begin{equation*}
	\{(0)\}, \quad \{ (E), \,(F) \},\quad \{ (C-\beta)\,|\,\beta\in \mK^* \},\quad \{ (E,F) \}, \quad \{ (E,F,K-\gamma)\,|\, \gamma\in \mK^* \}. 
	\end{equation*}
\end{corollary}
\begin{proof}
	The result follows from Theorem \ref{B16Oct21} and Theorem \ref{A24Oct21} immediately. 
\end{proof}

In the next theorem, we calculate the group of automorphisms  of the algebra $D_q$. 

\begin{theorem} \label{25Oct21} 
	$\Aut\,D_q \simeq (\mK^*)^4\rtimes {\rm SL}_2(\Z)$.  More precisely, any $\rho \in \Aut \,D_q$ is of the following form
	\begin{equation*}
	\begin{aligned}
	\rho_{\l, \mu, \gamma,\nu,i,j,m,n}: \quad &K \mapsto \l K^{i}a^{j} , \,\,\,\, a \mapsto \mu K^m a^{n},  \,\,\,\, F \mapsto \gamma K^{-i+2m+1}a^{-j+2n-2} F,\\
	& b \mapsto \mu \gamma^{-1} q^{ij+2mn-3in-2j+3n}K^{i-m-1} a^{j-n+1}b,\\
	& E \mapsto \nu K^{-2m}a^{-2n+2} E, \,\,\,\,\,\,\,c \mapsto \l \mu^{-1} \nu^{-1}q^{in+3mn-j-n}K^{i+m-1}a^{j+n-1}c
	\end{aligned}
	\end{equation*}
	where $(\l, \mu,\gamma, \nu) \in (\mK^*)^4$ and $\begin{psmallmatrix} i & j \\ m &n 
	\end{psmallmatrix}\in {\rm SL}_2(\Z)$. 
\end{theorem}

\begin{proof}
	It is straightforward to verify that the map $\rho_{\l, \mu, \gamma,\nu,i,j,m,n}$ defines an automorphism of $D_q$. The subgroup $G$ of $\Aut\,D_q$ generated by $\rho_{\l, \mu, \gamma,\nu,i,j,m,n}$ ($(\l, \mu, \gamma, \nu) \in (\mK^*)^4$ and $\begin{psmallmatrix} i & j \\ m &n 
	\end{psmallmatrix}\in {\rm SL}_2(\Z)$) is isomorphic to $(\mK^*)^4 \rtimes {\rm SL}_2(\Z)$. It remains to show that any automorphism of the algebra $D_q$ is of the form $\rho_{\l, \mu, \gamma,\nu,i,j,m,n}$. 
	First, we note that the group of units of $D_q$ is equal to $\{ \mK^* K^i a^j \,|\,i,j\in\Z \}$. Since any automorphism $\rho$ of $D_q$ preserves the group of units, we have $\rho(K)=\l K^i a^j$ and $\rho(a)=\mu K^m a^n$ for some $\l, \mu \in \mK^*$ and $i,j,m,n \in \Z$. By identification of coefficients in the identity $\rho(K)\rho(a)=q^{-1}\rho(a) \rho(K)$ we deduce that $q^{in-jm-1}=1$, thus $in-jm=1$ since $q$ is not a root of unity.  By Theorem \ref{A19Oct21}, the height one primes of $D_q$ are $(\phi)$ and $(\psi)$. So there are two options, either $(\phi)$ and $(\psi)$ are $\rho$-invariant ideals, or otherwise, they are interchanged under the action of $\rho$. We shall see that the second case is impossible. 
	
	Assume first that the ideals
	$(\phi)$ and $(\psi)$ are $\rho$-invariant.  Since $\phi$ and $\psi$ are normal elements,  $\rho(\phi)=\alpha K^s a^t \phi$ and $\rho(\psi)=\alpha' K^{u}a^{v}\psi$ for some $\alpha, \alpha'\in \mK^*$ and $s, t, u, v\in \Z$. The following calculation determines $s,t$, $u,v$ (in terms of $i,j,m,n$). In fact,  by identification of coefficients in the equalities  $\rho(\phi) \rho(K)=q \rho(K) \rho(\phi)$ and  $\rho(\phi) \rho(a)=\rho(a)\rho(\phi)$, we deduce that
	$js-it=i-1$ and $ns-mt=m$, respectively. Since $in-jm=1$, solving this system gives 
	$s=m$ and $t=n-1$. Thus $\rho(\phi)=\alpha K^{m} a^{n-1} \phi$. In a similar manner, by identification of coefficients in the equalities $\rho(\psi) \rho(K)=q^{-1} \rho(K) \rho(\psi)$ and $\rho(\psi) \rho(a)=q^{-1}\rho(a) \rho(\psi)$, we deduce that $ju-iv=1-i-j$ and $nu-mv=1-m-n$. 
	Solving this system gives $u=i-m-1$ and $v=j-n+1$, thus $\rho(\psi)=\alpha' K^{i-m-1}a^{j-n+1}\psi$. 
	Let us now consider the actions of $\rho$ on the generators $F$ and $b$. 
	From the expression of $\phi$, we have $\rho(\phi)=(1-q^2)\rho(F)\rho(b)+q^2 \rho(a)$. Thus 
	\begin{equation*}
	(1-q^2)\rho(F)\rho(b)=\rho(\phi)-q^2\rho(a)=(1-q^2)\alpha K^m a^{n-1}Fb+q^2(\alpha-\mu)K^ma^n. 
	\end{equation*}
	Since $D_q$ is a domain, for any $x\in \{F,E,b,c\}$ and $u,v\in D_q$,  we have $\deg_x(uv)=\deg_x(u) +\deg_x(v)$. Then from the above equality we deduce that $\rho(F)=\sum \gamma_{\ell_1,\ell_2, \epsilon_1, \epsilon_2} K^{\ell_1}a^{\ell_2}F^{\epsilon_1}b^{\epsilon_2}$ where $\ell_1, \ell_2 \in\Z$, $\gamma_{\ell_1,\ell_2, \epsilon_1, \epsilon_2}\in \mK$, $\epsilon_1, \epsilon_2\in \N$,  and $\epsilon_1 \leq 1$, $\epsilon_2 \leq 1$. By identification of coefficients of the identities $\rho(a)\rho(F)=q^{-1} \rho(F) \rho(a)$ and $\rho(\phi)\rho(F)=q \rho(F) \rho(\phi)$, we obtain, respectively, that
	\begin{subequations} 
		\begin{align}
		&n \ell_1-m \ell_2-(n+2m)\epsilon_1+(m+n) \epsilon_2 =-1, \label{Eq1} \\
		&n \ell_1 -m \ell_2-(n+2m-2) \epsilon_1+(n+m-2) \epsilon_2=1. \label{Eq2} 
		\end{align}
	\end{subequations}
	Subtracting (\ref{Eq1}) from (\ref{Eq2}) yields that  $\epsilon_1-\epsilon_2=1$. As $\epsilon_1 \leq 1$ and $\epsilon_2 \leq 1$, we obtain $\epsilon_1=1$, $\epsilon_2=0$. It follows that $\rho(F)=\sum \gamma_{\ell_1, \ell_2} K^{\ell_1}a^{\ell_2}F$, and $\rho(b)=\sum \gamma_{e_1, e_2} K^{e_1}a^{e_2}b$ where $e_1, e_2\in \Z$.  The equalities (\ref{Eq1}) and (\ref{Eq2}) now become
	\begin{equation} \label{Eq3} 
	n \ell_1-m \ell_2=n+2m-1. 
	\end{equation}
	By identification of coefficients in the equality $\rho(K) \rho(F)=q^{-2} \rho(F) \rho(K)$, we obtain
	\begin{equation}\label{Eq4} 
	j \ell_1-i \ell_2=2i +j-2. 
	\end{equation}
	Combining the equalities (\ref{Eq3}) and (\ref{Eq4}), we deduce that $\ell_1=-i+2m+1$ and $\ell_2=-j+2n-2$. Thus $\rho(F)=\gamma K^{-i+2m+1} a^{-j+2n-2}F$ for some $\gamma \in \mK^*$. Similarly, by identification of coefficients in the equalities $\rho(K)\rho(b)=q \rho(b) \rho(K)$ and $\rho(a) \rho(b)=q\rho(b) \rho(a)$ we obtain, respectively, that
	\begin{equation*}
	j e_1-i e_2=-j-i+1, \quad n e_1-m e_2=-n-m+1. 
	\end{equation*} 
	Thus $e_1=i-m-1$ and $e_2=j-n+1$, and $\rho(b)$ can be written as $\rho(b)=\gamma' K^{i-m-1}a^{j-n+1} b$ for some $\gamma' \in \mK^*$. Now substituting the expressions of $\rho(F)$, $\rho(a)$ and $\rho(b)$ into the identity $\rho(F)\rho(b)=q^{-1}\rho(b) \rho(F)+\rho(a)$, by identification of coefficients we deduce that $\gamma'=\mu \gamma^{-1}q^{ij+2mn-3in-2j+3n}$. 
	Let us consider the actions of $\rho$ on the generators $c$ and $E$. 
	Since the ideals $(\phi)$ and $(\psi)$ are $\rho$-invariant, from the diagram in Theorem \ref{A19Oct21}, we see that the automorphism $\rho$ preserves the set $\{ \gm_{\beta}=(\phi-\beta cK^{-1}, \psi)\,|\, \beta\in \mK^* \}$ of maximal ideals of $D_q$ over $(\psi)$. Thus $\rho(\gm_{\beta})=\gm_{\beta'}$ for some $\beta' \in \mK^*$.  Combining this with the identity $\rho(\psi)=(1-q^{-2})\rho(E) \rho(c)+q^{-2}\rho(a)^{-1} \rho(K)$, we get $\rho(c)=\delta K^{i+m-1}a^{j+n-1}c$ for some $\d \in \mK^*$, and so $\rho(E)$ is of the form $\rho(E)=\sum \nu_{r_1, r_2} K^{r_1}a^{r_2}E$ for some $r_1,r_2\in \Z$ and $\nu_{r_1, r_2}\in \mK$. Then identification of coefficients in the equalities $\rho(K)\rho(E)=q^2\rho(E)\rho(K)$ and $\rho(a)\rho(E)=\rho(E)\rho(a)$ yields, respectively, that $jr_1-ir_2+2i=2$ and $nr_1-mr_2+2m=0$. Therefore, $r_1=-2m$, $r_2=-2n+2$, and consequently, $\rho(E)=\nu K^{-2m}a^{-2n+2}E$ for some $\nu \in \mK^*$. 
	 Finally,  from the identity $\rho(E)\rho(c)=\rho(c)\rho(E)+\rho(a)^{-1}\rho(K)$ we deduce that $\nu \d=\l \mu^{-1}q^{in+3mn-j-n}$. Therefore, $\rho(c)=\l \mu^{-1} \nu^{-1}q^{in+3mn-j-n}K^{i+m-1}a^{j+n-1}c$.  This proves that $\rho$ is of the form $\rho_{\l, \mu, \gamma,\nu,i,j,m,n}$ as given in the theorem.

	Suppose now that the ideals $(\phi)$ and $(\psi)$ are interchanged under an automorphism $\tau$ of $D_q$, we seek a contradiction. Recall that $\tau(K)=\l K^i a^j$ and $\tau(a)=\mu K^m a^n$ where $in-jm=1$. Since the elements $\phi$ and $\psi$ are normal in $D_q$, we have $\tau(\phi)=\nu K^{s_1}a^{s_2}\psi$ for some $\nu \in \mK^*$ and $s_1,s_2\in \Z$. By identification of coefficients in the equalities $\tau(K) \tau(\phi)=q^{-1} \tau(\phi) \tau(K)$ and $\tau(a)\tau(\phi)=\tau(\phi) \tau(a)$, we obtain $js_1-i s_2=-i-j-1$ and $ns_1-ms_2=-m-n$, respectively. Solving this system gives $s_1=m-1$ and $s_2=n+1$, hence $\tau(\phi)=\nu K^{m-1}a^{n+1}\psi$. Substituting the explicit expression of $\phi$ and $\psi$ into this identity yields that
	\begin{equation*}
	(1-q^2)\tau(F) \tau(b)=\nu (1-q^{-2}) K^{m-1}a^{n+1}Ec+(\nu q^{n-2}-\mu q^2) K^m a^n. 
	\end{equation*}
	Comparing the degrees of the generators $E$, $c$, $F$, $b$ on both sides, we see that $\tau(F)$ is of the form $\tau(F)=\sum \gamma_{e_1, e_2,\varepsilon_1, \varepsilon_2}K^{e_1}a^{e_2}E^{\varepsilon_1}c^{\varepsilon_2}$ where $e_1, e_2\in \Z$, $\varepsilon_1, \varepsilon_2 \in \N$ and $\varepsilon_1, \varepsilon_2 \leq 1$. By identification of coefficients in the equalities $\tau(a) \tau(F)=q^{-1} \tau(F) \tau(a)$ and $\tau(\phi) \tau(F)=q \tau(F) \tau(\phi)$, we obtain
	\begin{equation*}
	\begin{aligned}
	&ne_1-me_2+2m \varepsilon_1+(n-m)\varepsilon_2 =-1, \\
	&n e_1-me_2+2(m-1)\varepsilon_1+(n-m+2) \varepsilon_2=1. 
	\end{aligned}
	\end{equation*}
	Subtracting the first equality from the second one yields that $\varepsilon_2-\varepsilon_1=1$. Thus $\varepsilon_1=0$ and $\varepsilon_2=1$, and so $ne_1-me_2+n-m=-1$. Now by identification of coefficients in the equality $\tau(K) \tau(F)=q^{-2} \tau(F) \tau(K)$ we obtain that $je_1-i e_2+j-i=-2$. Solving for $e_1$ and $e_2$ we obtain $e_1=2m-i-1$ and $e_2=2n-j-1$. Thus $\tau(F)=\gamma K^{2m-i-1}a^{2n-j-1}c$ for some $\gamma \in \mK^*$. Now $\tau(b)$ is of the form $\tau(b)=\sum \varsigma_{\ell_1, \ell_2} K^{\ell_1} a^{\ell_2}E$ where $\ell_1, \ell_2 \in \Z$. Identification of coefficients in the equalities $\tau(K) \tau(b)=q \tau(b) \tau(K)$ and $\tau(a) \tau(b)=q \tau(b) \tau(a)$ yields, respectively, that $j \ell_1 - i \ell_2=-2i+1$ and $n \ell_1- m \ell_2 =-2m +1$. Solving this system gives $\ell_1=i-m$ and $\ell_2=j-n+2$. Thus $\tau(b)=\varsigma K^{i-m}a^{j-n+2}E$ for some $\varsigma \in \mK^*$. To summarize, now we have 
	\begin{equation*}
	\begin{aligned}
	\tau(K)&=\l K^i a^j, &\tau(a)&=\mu K^m a^n, \\
	\tau(F)&=\gamma K^{2m-i-1}a^{2n-j-1}c, &\tau(b)&=\varsigma K^{i-m}a^{j-n+2}E. 
	\end{aligned}
	\end{equation*}
	Note that $\tau^2$ is an automorphism of $D_q$ that preserves the ideals $(\phi)$ and $(\psi)$. From the expression of $\tau^2(b)$ and the result obtained in the previous case, one sees that $\tau(E)=\d K^{t_1}a^{t_2}b$ for some $\d \in \mK^*$ and $t_1, t_2 \in \Z$. Then identification of coefficients in the equalities $\tau(a) \tau(E)=\tau(E)\tau(a)$ and $\tau(K)\tau(E)=q^2\tau(E) \tau(K)$ yields, respectively, that 
	\begin{equation} \label{mnij} 
	n t_1-mt_2+m+n=0, \quad j t_1 -i t_2 +i+j=2. 
	\end{equation}
	Furthermore, by identification of coefficients in the identity $\tau(F)\tau(E)=\tau(E)\tau(F)$, one deduces that $(2n-j)t_1-(2m-i)t_2+2(m+n)-(i+j)=2$. Substituting (\ref{mnij}) into this equality one obtains $0=4$, a contradiction. 
\end{proof}

\begin{corollary} \label{a11Nov21} 
	The action of $\Aut\,D_q$ on $\Spec (D_q)$ has exactly 6 orbits:
	\begin{equation*}
	\{ (0) \}, \quad \{ (\phi) \}, \quad \{(\psi) \}, \quad \{(\phi,\psi) \}, \quad \{ (\phi, \psi-\alpha b)\,|\, \alpha \in \mK^* \}, \quad \{ (\phi-\beta cK^{-1}, \psi)\,|\, \beta \in \mK^* \}. 
	\end{equation*}	
\end{corollary}
\begin{proof}
	The result follows from Theorem \ref{A19Oct21} and Theorem \ref{25Oct21}. 
\end{proof}

\section{Simple weight modules of $D_q(E_2)$} \label{SimpleWeight} 

The aim of this section is to give a classification of simple weight $D_q$-modules. The structure and representations of the centralizer $C(K,D_q)$ of the element $K$ in $D_q$  are  key ingredients in this procedure. The problem of classification of simple weight $D_q$-modules is essentially reduced to the problem of classification of all $C(K,D_q)$-modules (Theorem \ref{D17Nov21}). All the primitive factors of $C(K,D_q)$ are explicitly described as generalized Weyl algebras (Proposition \ref{d16Nov21}), and using this fact, a classification of simple unfaithful $C(K,D_q)$-modules is obtained.

The following proposition gives an explicit description of the centralizer $C(K, D_q)$ of $K$ in $D_q$ via generators and defining relations. 

\begin{proposition} \label{a16Nov21} 
	Let	$x_1:=a^2E$, $y_1:=a^{-1}c$, $x_2:=a^{-2}F$ and $y_2:=ab$.
	\begin{enumerate}
		\item  The algebra $\sC:=C(K, D_q)$ is  generated by the elements $K^{\pm 1}$, $x_1$, $y_1$, $x_2$ and $y_2$  subject to the following defining relations:	
		\begin{equation} \label{Centralizer} 
		\begin{aligned}
		&x_1y_1-q^2 y_1 x_1=K, &\quad  &x_2y_2-q^{-2}y_2 x_2=q, \\
		&x_2 x_1=q^2 x_1 x_2, & &y_2 x_1=q^{-2}x_1 y_2, \\
		&x_2 y_1=q^{-2}y_1 x_2, &&y_2 y_1=q^2 y_1y_2, \\
		& K K^{-1}=1=K^{-1}K, && [K,\,\, \cdot]=0. 
		\end{aligned}
		\end{equation}
		In particular, $\sC$ is a  Noetherian domain with $\GK(\sC)=5$. 
		\item The algebra $D_q$ is free as a left and right $\sC$-module, $D_q=\bigoplus_{i\in \Z} \sC a^i = \bigoplus_{i\in \Z} a^i \sC. $
	\end{enumerate}	
\end{proposition}
\begin{proof}
	The algebra $D_q$ admits a PBW type basis, and the elements $K^i a^j E^{k_1}c^{k_2}F^{k_3}b^{k_4}$ ($i,j\in \Z$, and $k_1,\ldots, k_4 \in \N$) form a basis of $D_q$.
	Note that $D_q$ is a $\Z$-graded algebra with respect to the adjoint action of $K$; that is, $D_q=\oplus_{n\in \Z} D_q[n]$ where $D_q[n]=\{ x\in D_q\,|\, KxK^{-1}=q^n x \}$. Clearly, the algebra  $\sC$ is exactly the degree zero component $D_q[0]$ with respect to this grading. Then
	\begin{equation*}
	\begin{aligned}
	\sC=&\bigoplus_{\stackrel{-j+2k_1-k_2-2k_3+k_4=0}{i,j\in \Z, k_1, \ldots k_4\in \N}} \mK K^i a^j E^{k_1}c^{k_2}F^{k_3}b^{k_4}\\
	&= \bigoplus_{i\in \Z, k_1, \ldots k_4\in \N} \mK K^i a^{2k_1-k_2-2k_3+k_4}E^{k_1}c^{k_2}F^{k_3}b^{k_4}\\
	&= \bigoplus_{i\in \Z, k_1, \ldots k_4\in \N} \mK K^i (a^2E)^{k_1}(a^{-1}c)^{k_2}(a^{-2}F)^{k_3}(ab)^{k_4}\\
	&=\bigoplus_{i\in \Z, k_1, \ldots k_4\in \N} \mK K^i x_1^{k_1}y_1^{k_2}x_2^{k_3}y_2^{k_4}. 
	\end{aligned}
	\end{equation*}
	Therefore $\sC$ is generated by $K^{\pm 1}$, $x_1, y_1, x_2, y_2$, and  the elements $K^i x_1^{k_1}y_1^{k_2}x_2^{k_3}y_2^{k_4}$ ($i\in \Z$, and $k_1, \ldots, k_4\in \N$) form a basis of $\sC$. 
	It is straightforward to verify that the relations in (\ref{Centralizer}) hold in the algebra $\sC$. Let $\mathscr{R}$ be the abstract $\mK$-algebra generated by $K^{\pm 1}$, $x_1$, $y_1$, $x_2$ and $y_2$ subject to the defining relations as in (\ref{Centralizer}). Then there is a natural epimorphism of algebras $\pi: \mathscr{R}\twoheadrightarrow \sC$. It is easy to see that $\mathscr{R}$ has a PBW type basis. Comparing the bases of $\mathscr{R}$ and $\sC$ yields that $\pi$ is an isomorphism. This proves that the relations in (\ref{Centralizer}) are defining relations of $\sC$. It is easy to see that the algebra $\sC$ can be presented as an iterated Ore extension of the form $\sC=\mK[K^{\pm 1},y_1][y_2, \s_1][x_1;\s_2, \d_2][x_2;\s_3,\d_3]$, and in particular $\sC$ is a Noetherian domain with $\GK(\sC)=5$. 
	
	2. Statement 2 follows from the bases of $D_q$ and $\sC$ given in the proof of statement 1. 	
\end{proof}

Set $u:=(q^2-1)y_1x_1+K$ and $v:=(q^{-3}-q^{-1})y_2x_2+1$. Then $u,v \in \sC$, and an easy computation shows that $u=a\psi$ and $v=a^{-1}\phi$. The elements $u$ and $v$ are normal in $\sC$, in fact
\begin{equation}\label{uvnorm} 
\begin{aligned}
x_1 u &=q^2 ux_1,  &x_2 u&=u x_2, &y_1 u&= q^{-2} u y_1, &y_2 u&=u y_2, \\
x_1 v&=v x_1,  &x_2 v&=q^{-2}v x_2, &y_1 v&=v y_1, &y_2 v&=q^2 vy_2. 
\end{aligned}
\end{equation} 
The following lemma shows that the centre of $\sC$ is the Laurent polynomial algebra $\mK[K, K^{-1}]$. 
\begin{lemma} \label{b16Nov21} 
	$Z(\sC)=\mK[K, K^{-1}]$. 
\end{lemma}
\begin{proof}
	Let $\sC X^{-1}$ be the localization of $\sC$ at the Ore set $X=\{\mK^*x_1^i x_2^j\,|\, i,j\in \N \}$. After a change of generators, namely $y_1 \rightsquigarrow u$ and $y_2 \rightsquigarrow v$, we see that $\sC X^{-1}$ can be presented as a tensor product of algebras
	\begin{equation} \label{CX} 
	\sC X^{-1}=\mK[K, K^{-1}] \otimes \mK_{\bf{q}}[x_1^{\pm 1}, x_2^{\pm 1}, u,v]
	\end{equation}	
	where the second tensor component is a quantum polynomial algebra and the $4\times 4$ matrix $\bf{q}$ can be  easily obtained from the commutation relations in (\ref{uvnorm}). It is easily deduced that the centre of the second tensor component is trivial, thus $Z(\sC X^{-1})=\mK[K, K^{-1}]$, and therefore $Z(\sC)=\sC \cap Z(\sC X^{-1})=\mK[K, K^{-1}]$. 
\end{proof}

For any $\chi \in \mK^*$, let us define
$\sC(\chi):=\sC/\sC(K-\chi). $
Then the algebra $\sC(\chi)$ is generated by the elements $x_1$, $y_1$, $x_2$, $y_2$ subject to the relations in (\ref{Centralizer}) by setting $K=\chi$. So $\sC(\chi)$ is a `mixture' of two quantum Weyl algebras. Clearly, $\sC(\chi)$ is a Noetherian domain with $\GK (\sC(\chi))=4$. 

\begin{lemma} \label{c16Nov21} 
	In the algebra $\sC(\chi)$, for any positive integer $i \geq 1$, we have 
	\begin{equation*}
	(x_1^i)=(y_1^i)=(x_2^i)=(y_2^i)=\sC(\chi). 
	\end{equation*}	
\end{lemma}
\begin{proof}
	We only show that $(x_1^i)=\sC(\chi)$, since other cases can be proved similarly. From the identity $x_1 y_1-q^2 y_1x_1=\chi$ it follows immediately that $(x_1)=\sC(\chi)$. Note that in $\sC(\chi)$ we have 
	\begin{equation*}
	x_1^i y_1-q^{2i}y_1 x_1^i=\chi \tfrac{1-q^{2i}}{1-q^2} x_1^{i-1}. 
	\end{equation*}
	The above equality yields that $(x_1^i)=(x_1^{i-1})=\cdots=(x_1)=\sC(\chi)$. 
\end{proof}

The following proposition describes certain factor algebras of $\sC(\chi)$ which in fact are all the prime factors of $\sC(\chi)$. All these factor algebras are presented as classical GWAs. 
\begin{proposition}\label{d16Nov21} 
	Let $\chi \in \mK^*$, $\theta:=vx_1$, and $\omega:=ux_2$. 
	\begin{enumerate}
		\item The factor algebra $\sC(\chi)/(u)$ is a domain that can be presented as a GWA of the form
		\begin{equation*}
		\sC(\chi)/(u) \simeq \mK[v, x_1^{\pm}][x_2, y_2;\s, \emph{\texttt{a}}=(q^{-3}-q^{-1})^{-1}(v-1)]
		\end{equation*}
		where $\s$ is the automorphism of $\mK[v,x_1^{\pm 1}]$ defined by $\s(v)=q^{-2}v$ and $\s(x_1)=q^2 x_1$. Moreover, the centre of $\sC(\chi)/(u)$ is the polynomial algebra $\mK[\theta]$. 
		
		\item The factor algebra $\sC(\chi)/(v)$ is a domain that can be presented as a GWA of the form
		\begin{equation*}
		\sC(\chi)/(v)\simeq \mK[u, x_2^{\pm 1}][x_1,y_1;\s, \emph{\texttt{a}}=(q^2-1)^{-1}(u-\chi)]
		\end{equation*}
		where $\s$ is the automorphism of $\mK[u,x_2^{\pm 1}]$ defined by $\s(u)=q^2 u$ and $\s(x_2)=q^{-2}x_2$. Moreover, the centre of $\sC(\chi)/(v)$ is the polynomial algebra $\mK[\omega]$. 
		
		\item For any $\alpha \in \mK$, the factor algebra $\sC(\chi)/(u,\theta-\alpha)$ is a simple domain that can be expressed as a GWA of the following form
		\begin{equation*}
		\sC(\chi)/(u,\theta-\alpha) \simeq \mK[x_1^{\pm 1}][x_2, y_2;\s, \emph{\texttt{a}}=(q^{-3}-q^{-1})^{-1}(\alpha x_1^{-1}-1)]
		\end{equation*}
		where $\s$ is the automorphism of $\mK[x_1^{\pm 1}]$ defined by $\s(x_1)=q^2 x_1$. In particular, if $\alpha =0$ then $\sC(\chi)/(u, \theta)$ is isomorphic to the quantum torus $\mK_{q^{-2}}[x_1^{\pm 1}, x_2^{\pm 1}]$. 
		
		\item For any  $\beta \in \mK$, the factor algebra $\sC(\chi)/(v, \omega-\beta)$ is a simple domain that can be presented as a GWA of the following form
		\begin{equation*}
		\sC(\chi)/(v,\omega-\beta) \simeq \mK[x_2^{\pm 1}][x_1, y_1;\s, \emph{\texttt{a}}=(q^2-1)^{-1}(\beta x_2^{-1}-\chi)]
		\end{equation*}
		where $\s$ is the automorphism of $\mK[x_2^{\pm 1}]$ defined by $\s(x_2)=q^{-2}x_2$. 
	\end{enumerate}	
\end{proposition}
\begin{proof}
	1. From the expression of the element $u=(q^2-1)y_1 x_1+\chi$, we see that the elements $x_1$ and $y_1$ are invertible in $\Gamma:=\sC(\chi)/(u)$, and that $y_1=\chi(1-q^2)^{-1}x_1^{-1}$. Then $\Gamma$ is generated by $x_1^{\pm 1}$, $x_2$ and $y_2$ which is obviously can be expressed as a GWA as the form given in the statement. In particular, $\Gamma$ is a Noetherian domain of GK dimension 3.  It is easily verified that $\theta$ is central in $\Gamma$. Note that the localization $\Gamma[x_2^{-1}]$ of the algebra $\Gamma$ at the powers of $x_2$ can be expressed as a tensor product of algebras $\Gamma[x_2^{-1}]=\mK[\theta] \otimes \mK_{q^{-2}}[x_1^{\pm 1}, x_2^{\pm 1}]$ and the second tensor component is a central simple quantum torus. It follows that $Z(\Gamma[x_2^{-1}])=\mK[\theta]$, and consequently, $Z(\Gamma)=\mK[\theta]$.  
	
	2. The proof is similar to that of statement 1. 
	
	3. The fact that $\Gamma(\chi,\alpha):=\sC(\chi)/(u,\theta-\alpha)$ can be presented as a GWA as the given form follows immediately from statement 1. In particular, $\Gamma(\chi,\alpha)$ is a Noetherian domain of GK dimension 2. Notice that the localization $\Gamma(\chi,\alpha)[x_2^{-1}]$ of the algebra $\Gamma(\chi, \alpha)$ at the powers of $x_2$ is isomorphic to the quantum torus $\mK_{q^{-2}}[x_1^{\pm 1}, x_2^{\pm 1}]$, which is a simple algebra. Hence $\Gamma(\chi,\alpha)$ is a simple algebra iff $(x_2^i)=\Gamma(\chi, \alpha)$ for all $i \geq 1$. But this follows from Lemma \ref{c16Nov21}.  
	
	4. The proof is similar to that of statement 3. 
\end{proof}

Note that $\theta=vx_1=a\phi E$ and $\omega=ux_2=a^{-1}\psi F$. The following theorem gives an explicit description of the prime, completely prime, maximal and primitive ideals of $\sC(\chi)$ for any $\chi \in \mK^*$.
\begin{theorem} \label{A17Nov21} 
	For any $\chi \in \mK^*$, the prime spectrum of $\sC(\chi)$ is given below,
	\begin{equation*}
	\Spec(\sC(\chi))=\{  (0), \,(u),\,(v),\,(u,v) \} \,\sqcup\,\{ (u,\theta-\alpha)\,|\, \alpha \in \mK^* \}\,\sqcup\,\{ (v,\omega-\beta)\,|\, \beta \in \mK^* \}. 
	\end{equation*}
	The inclusions of prime ideas of $\sC(\chi)$ are depicted in the following diagram
	\begin{equation*}
	\begin{tikzpicture}[text centered]
	\node (vomega) at (2,2) {\qquad\qquad $\{(v, \,\omega-\beta)|\beta \in \mK^*\}$};
	\node (utheta) at (-2,2) {$\{(u,\,\theta-\alpha)|\alpha \in \mK^*\}$\qquad\qquad\qquad};
	\node (uv) at (0,2) {$(u,v)$};
	\node (v) at (1,1)  {$(v)$};
	\node (u)  at (-1,1)  {$(u)$};
	\node (0)  at (0,0)  {$(0)$};
	\draw [shorten <=-2pt, shorten >=-2pt] (vomega) -- (v) -- (0) -- (u)--(utheta);
	\draw [shorten <=-2pt, shorten >=-2pt] (u)--(uv)--(v);
	\end{tikzpicture} 
	\end{equation*} 	
	Furthermore, 
	\begin{enumerate}
		\item All prime ideals of $\sC(\chi)$ are completely prime. 
		\item $\Max(\sC(\chi))=\{ (u,\,\theta-\alpha)\,|\,\alpha \in \mK^* \} \, \sqcup\,\{(u,v)\}\,\sqcup\,\{ (v,\,\omega-\beta)\,|\, \beta\in \mK^* \}$. 
		\item $\Prim(\sC(\chi))=\{(0)\}\,\sqcup\,\Max(\sC(\chi))$. 
	\end{enumerate}
\end{theorem}

\begin{proof}
	Let $\CR:=\sC(\chi)X^{-1}$ be the localization of $\sC(\chi)$ at the Ore set $X=\{ \mK^*x_1^i x_2^j\,|\,i,j\in \N \}$. It follows from Lemma \ref{c16Nov21} that any prime ideal of $\sC(\chi)$ has no intersection with $X$, and thus the prime spectrum of $\sC(\chi)$ is homeomorphic to the prime spectrum of $\CR$. Note that the elements $u$ and $v$ are normal in $\CR$.  If we denote by $\CR[u^{-1}, v^{-1}]$ the localization of $\CR$ at the Ore set $\{ u^iv^j \,|\,i,j\in \N \}$, then by (\ref{CX}), $\CR[u^{-1},v^{-1}]=\mK_{\bf q}[x_1^{\pm}, x_2^{\pm},u^{\pm 1},v^{\pm 1} ]$ is a simple quantum torus. Consequently, for any non-zero prime ideal $\gp$ of $\CR$, either $u\in \gp$ or $v\in \gp$. In other words, the set of non-zero prime ideals of $\CR$ decomposes into a union of the sets $T_u=\{ \gp \in \Spec(\CR)\,|\,u\in \gp \}$ and $T_v=\{ \gp \in \Spec(\CR)\,|\, v\in \gp \}$. 
	For a subset $\CG$ of $\CR$, let us denote by $(\CG)_X$ the two-sided ideals of $\CR$ generated by $\CG$.
	The prime ideals of $\CR$ that contain $u$ are in bijection with the prime ideals of $\CR/(u)_X$. From Proposition \ref{d16Nov21}.(1), we see that $\CR/(u)_X \simeq [\sC(\chi)/(u)]X^{-1}\simeq \mK[\theta] \otimes \mK_{q^{-2}}[x_1^{\pm 1}, x_2^{\pm 1}]$. Since the quantum torus $\mK_{q^{-2}}[x_1^{\pm 1}, x_2^{\pm 1}]$ is a simple algebra, the prime spectrum of $\CR/(u)_X$ is homeomorphic to the prime spectrum of $\mK[\theta]$. Therefore $T_u=\{(u)_X, (u,\theta-\alpha)_X\,|\,\alpha \in \mK \}$. A similar argument shows that $T_v=\{ (v)_X,(v, \omega-\beta)_X\,|\, \beta\in \mK \}$. 
	
	Clearly $(0)_X \cap \sC(\chi)=(0)$. Let us show that $(u)_X \cap \sC(\chi)=(u)$. The inclusion ``$\supseteq$" is obvious. For the converse, if $z\in (u)_X \cap \sC(\chi)$ then $z x_1^i x_2^j\in (u)$ for some $i,j\in \N$. By Proposition \ref{d16Nov21}.(1), $(u)$ is a completely prime ideal of $\sC(\chi)$ and $x_1, x_2 \notin (u)$ thus $z\in (u)$.	Similarly, using Proposition \ref{d16Nov21}.(2) one obtains $(v)_X \cap \sC(\chi)=(v)$. For any $\alpha, \beta\in \mK$, we have $(u,\theta-\alpha)_X \cap \sC(\chi)=(u,\theta-\alpha)$ and $(v,\omega-\beta)_X \cap \sC(\chi)=(v, \omega-\beta)$, since the ideals $(u,\theta-\alpha)$ and $(v, \omega-\beta)$ are maximal in $\sC(\chi)$, see Proposition \ref{d16Nov21}.(3) and (4). It is easy to see that $(u, \theta)=(u,v)=(v, \omega)$.  All  prime ideals of $\sC(\chi)$ are shown in the diagram and the containments are obvious.

	1. This follows from Proposition \ref{d16Nov21}, since all prime factors of $\sC(\chi)$	are domains. 
	
	2. This follows from the diagram in the statement. 
	
	3. All maximal ideals of $\sC(\chi)$ are primitive. The zero ideal $(0)$ is primitive since it is locally closed. The ideals $(u)$ and $(v)$ are not primitive, since the factor algebras $\sC(\chi)/(u)$ and $\sC(\chi)/(v)$ have non-trivial centres, see Proposition \ref{d16Nov21}.(1) and (2).	
\end{proof}

\begin{corollary} \label{c17Nov21} 
	The algebra $\sC$ has no finite-dimensional modules. 
\end{corollary}
\begin{proof}
	If $\sC$ has a finite-dimensional module, then $\sC$ has a finite-dimensional simple module $M$. By Schur's Lemma, the central element $K$ of $\sC$ acts on $M$ by a non-zero scalar $\chi$, since the field $\mK$ is algebraically closed. Then $M$ is a simple $\sC(\chi)$-module, and the annihilator $\ann_{\sC(\chi)}(M)$ of $M$ is a primitive ideal of $\sC(\chi)$ with finite codimension. However, from Proposition \ref{d16Nov21}, all primitive factors of $\sC(\chi)$ are infinite-dimensional, a contradiction. 
\end{proof}

Now we  apply the results obtained for the centralizer $\sC$ to investigate weight modules over the algebra $D_q$. 

\begin{definition}
	A $D_q$-module $M$ is called a \emph{weight module} provided that $M$ decomposes into a direct sum of weight spaces; that is $M=\oplus_{\mu\in \mK^*}M_{\mu}$ where $M_{\mu}=\{x\in M\,|\,Kx=\mu x \}$. We denote by $\supp(M):=\{\mu\in\mK^*\,|\, M_{\mu} \neq 0 \}$ the set of all weights of $M$. A weight module is called a Harish-Chandra module if all its weight spaces are finite-dimensional.
\end{definition}

\begin{theorem}\label{B17Nov21} 
	Let $M$ be a weight $D_q$-module. 
	\begin{enumerate}
		\item For any $\mu \in \supp(M)$, we have $\supp(M) \supseteq \{ q^i \mu\,|\,i\in\Z \}$. If, in addition, the module $M$ is simple then $\supp(M) = \{ q^i \mu\,|\,i\in\Z \}$.
		\item Every weight space of $M$ is infinite-dimensional. In particular, the algebra $D_q$ has no Harish-Chandra modules. 
	\end{enumerate}	
\end{theorem}
\begin{proof}
	1. Since the element $a$ is invertible in $D_q$,  the map $a^i: M_{\mu} \rightarrow M_{q^{-i}\mu}$ is a bijection for any $i\in \Z$ and $\mu \in \supp(M)$. This implies that  $\supp(M) \supseteq \{ q^i \mu\,|\,i\in\Z \}$ for any $\mu \in \supp(M)$. If $M$ is simple, then any non-zero weight vector $v\in M_{\mu}$ generates the module $M$, i.e. $M=D_q v$. Using the PBW basis of $D_q$ we see that $\supp(M) \subseteq \{ q^i \mu\,|\, i\in \Z \}$. Hence $\supp(M) = \{ q^i \mu\,|\,i\in\Z \}$.
	
	2. Notice that every weight space $M_{\mu}$ of $M$ is a module over the centralizer $\sC=C(K, D_q)$. The result then follows from Corollary \ref{c17Nov21}.  
\end{proof}

Let $R$ be an algebra and $M$ a left $R$-module. For an automorphism $\tau$ of $R$, we define a new left $R$-module $M^{\tau}$ as follows. As a vector space $M^{\tau}=M$, the action of $R$ on $M^{\tau}$ is defined by $r\cdot m=\tau(r)m$ for $r\in R$ and $m\in M$. The module $M^{\tau}$ is called the twisting of $M$ by $\tau$. The normal element $a$ of $D_q$ induces an automorphism of $\sC$, namely $\tau_a: \sC \rightarrow \sC$, $t \mapsto a^{-1}t a$. On the generators of $\sC$, the automorphism $\tau_a$ is given by
\begin{equation} \label{Autau} 
\tau_a: \quad K \mapsto q^{-1}K, \quad x_1 \mapsto x_1, \quad y_1 \mapsto q^{-1}y_1, \quad x_2 \mapsto qx_2, \quad y_2 \mapsto q^{-1}y_2. 
\end{equation}

The following theorem provides a strategy to construct simple weight $D_q$-modules. 
It shows that the problem of classification of simple weight $D_q$-modules is essentially reduced to a problem of classification of all simple modules over the centralizer $\sC$ of the element $K$ in $D_q$.

\begin{theorem} \label{D17Nov21} 
	The following map defines a surjection 
	\begin{equation*}
	\widehat{\sC} \twoheadrightarrow \widehat{D_q}({\rm weight}), \quad M \mapsto \ind_{\sC}^{D_q}(M):=D_q \otimes_{\sC} M. 
	\end{equation*}
	Moreover, for $M, M' \in \widehat{\sC}$, $\ind_{\sC}^{D_q}(M) \simeq \ind_{\sC}^{D_q}(M')$ iff $_{\sC}M' \simeq {_{\sC}M^{\tau_a^i}}$ for some $i \in \Z$ where $\tau_a$ is the automorphism of $\sC$ defined in (\ref{Autau}). 
\end{theorem}
\begin{proof}
	Recall from Proposition \ref{a16Nov21}.(2) that $D_q$ is free as a right $\sC$-module, it follows that
	\begin{equation*}
	\ind_{\sC}^{D_q}(M)=\bigoplus_{i\in \Z} a^i \otimes M.
	\end{equation*}
	Since $M$ is a simple $\sC$-module, Schur's Lemma (the infinite-dimensional version due to Dixmier \cite{Dixmier-nilpotent}) implies that $(K-\chi)M=0$ for some $\chi \in \mK^*$. Then it is clear that $\ind_{\sC}^{D_q}(M)$ is a weight $D_q$-module whose support is equal to $\{ q^i \chi\,|\, i\in \Z \}$. Let $N$ be a non-zero submodule of $\ind_{\sC}^{D_q}(M)$, and pick a non-zero element $x=\sum a^i \otimes m_i \in N$ with each $m_i$ a non-zero element of $M$. Notice that the summands $a^i \otimes m_i$ of $x$ have distinct weights.  Thus each of them belongs to $N$, and hence $1\otimes m_i\in N$ since $a$ is invertible. The simplicity of $M$ as a $\sC$-module yields that $1\otimes M \subset N$ and hence $\ind_{\sC}^{D_q}(M)=N$, since $1\otimes M$ generates $\ind_{\sC}^{D_q}(M)$. This proves the simplicity of the $D_q$-module $\ind_{\sC}^{D_q}(M)$. Let us show that the map defined in the statement is a surjection. If $V$ is a simple weight $D_q$-module, then every weight space  $V_{\chi}$ of $V$ is a simple $\sC$-module. From the above argument we know that $\ind_{\sC}^{D_q}(V_{\chi})$ is a simple weight $D_q$-module. Since $\ind_{\sC}^{D_q}(V_{\chi}) \subseteq V$, the simplicity of $V$ then yields that $\ind_{\sC}^{D_q}(V_{\chi}) = V$. This proves that the map is a surjection.
	
	If $\ind_{\sC}^{D_q}(M) \simeq \ind_{\sC}^{D_q}(M')$ then, as a $\sC$-module, the weight space $M'$ is isomorphic to $a^i \otimes M$ for some $i \in \Z$. But the $\sC$-module $a^i \otimes M$ is clearly isomorphic to $_{\sC}M^{\tau_a^i}$, thus $_{\sC}M' \simeq {_{\sC}M^{\tau_a^i}}$. Conversely, if $_{\sC}M' \simeq {_{\sC}M^{\tau_a^i}}$ for some $i\in \Z$, then
	\begin{equation*}
	\ind_{\sC}^{D_q}(M') \simeq \bigoplus_{j\in \Z}a^j \otimes M^{\tau_a^i} \simeq \bigoplus_{j\in \Z} a^{j+i}\otimes M \simeq  \ind_{\sC}^{D_q}(M).
	\end{equation*}
	This completes the proof.
\end{proof}

According to Theorem \ref{D17Nov21}, the key ingredient in the construction of simple weight $D_q$-modules is the classification of simple $\sC$-modules. We now turn to this aspect. By Schur's Lemma,
\begin{equation*}
\widehat{\sC}=\bigsqcup_{\chi\in \mK^*} \widehat{\sC(\chi)}. 
\end{equation*}
Recall that an $R$-module $M$ is said to be faithful if its annihilator is zero. 
By Proposition \ref{d16Nov21}.(3) and (4), for any non-zero primitive ideal $\gp$ of $\sC(\chi)$, the factor algebra $\sC(\chi)/\gp$ is a simple domain, and therefore all $\sC(\chi)/\gp$-modules are faithful. In other words, the $\sC(\chi)/\gp$-modules are exactly the $\sC(\chi)$-modules with annihilators equal to $\gp$. It follows that 
\begin{equation} \label{Cchi} 
\widehat{\sC(\chi)} = \widehat{\sC(\chi)}({\rm unfaithful}) \;\sqcup \; \widehat{\sC(\chi)}({\rm faithful})
=\bigsqcup_{\gp \in \Prim \sC(\chi)\setminus\{(0)\}} \widehat{\sC(\chi)/\gp}\,\, \sqcup\,\, \widehat{\sC(\chi)}({\rm faithful}).
\end{equation}
By the explicit description of the primitive factors of $\sC(\chi)$ obtained in Proposition \ref{d16Nov21}, we are able to give a classification of simple unfaithful $\sC(\chi)$-modules. We shall also give several examples of simple faithful $\sC(\chi)$-modules. 

From Proposition \ref{d16Nov21}.(3) and (4), the primitive factors of $\sC(\chi)$ are quantum GWAs. In order to deal with the two cases simultaneously, we consider the GWAs of the following special form. Let $A:=A(\texttt{a}(h), q)$ be the ring generated by $\mK[h,h^{-1}]$, $x$ and $y$ subject to the relations
\begin{equation} \label{QGWA} 
xh=qhx, \,\, yh=q^{-1}h y, \,\, yx=\texttt{a}(h), \,\,xy=\texttt{a}(qh)
\end{equation}
where the Laurent polynomial $\texttt{a}(h) \in \mK[h,h^{-1}]$ is either a unit or has only one root, say $\zeta$, in $\mK^*$. In other words, either $\texttt{a}(h)=\mu h^i$ or $\texttt{a}(h)=\mu h^i (h-\zeta)$ for some $\mu\in \mK^*$ and $i\in\Z$. With this assumption, it is easily deduced that  $A(\texttt{a}(h),q)$ is a simple algebra.  The algebra $A$ is $\Z$-graded; that is, $A=\oplus_{i\in \Z}\mK[h^{\pm 1}]v_i$ where $v_i=x^i$ if $i\geq 0$ and $v_i=y^{-i}$ if $i<0$. Note that if $\texttt{a}(h)$ is a unit of $\mK[h,h^{-1}]$ then $A$ is isomorphic to the quantum torus $\mK_q[x^{\pm 1}, h^{\pm 1}]$. We denote by $\mathscr{O}(\zeta):=\{q^i \zeta\,|\, i\in \Z \}$ the orbit of $\zeta$ under the action of $q^{\Z}$ on $\mK^*$. 
By Proposition \ref{d16Nov21}.(3) and (4), 
\begin{equation} \label{Cchip} 
\begin{aligned}
{\rm if}\,\, \gp&=(u, \theta-\alpha), \,\,{\rm then}\,\,\sC(\chi)/\gp \simeq A\big((q^{-1}-q^{-3})^{-1}h^{-1}(h-\alpha), q^2 \big), \quad{\rm and}\\
{\rm if}\,\,\gp&=(v, \omega-\beta), \,\,{\rm then}\,\, \sC(\chi)/\gp \simeq A\big((1-q^2)^{-1}\chi h^{-1}(h-\chi^{-1}\beta), q^{-2} \big). 
\end{aligned}
\end{equation}


Let $R$ be an algebra, $M$ an $R$-module, and $S \subset R$ a left Ore set. The submodule $\tor_S(M):=\{ m\in M\,|\,sm=0 \,\,\text{for some}\,\, s\in S \}$ is called the $S$-torsion submodule of $M$. The module $M$ is said to be $S$-torsion if $\tor_S(M)=M$, and $S$-torsionfree if $\tor_S(M)=0$. If $S=\mK[h^{\pm 1}]\setminus \{0\}$ for some $h\in R$ then, by abusing the language, we simply say that $S$-torsion (resp. torsionfree) is $\mK[h^{\pm 1}]$-torsion (resp. torsionfree). Note that each simple $R$-module is either $S$-torsion or $S$-torsionfree. In particular, for the algebra $A=A(\texttt{a}(h),q)$ defined in (\ref{QGWA}), the set $\widehat{A}$ decomposes into a disjoint union of the following two subsets
\begin{equation} \label{DAmod} 
\widehat{A}=\widehat{A}\,(\mK[h^{\pm 1}]\text{-torsion}) \,\,\sqcup\,\,\widehat{A}\,(\mK[h^{\pm 1}]\text{-torsionfree}). 
\end{equation}

The following proposition gives an explicit description of the set $\widehat{A}\,(\mK[h^{\pm 1}]\text{-torsion})$.
\begin{proposition} \label{a19Nov21} 
	Let $M$ be a simple $\mK[h^{\pm1}]$-torsion module over the algebra $A=A(\emph{\texttt{a}}(h),q)$.
	\begin{enumerate}
		\item If $\emph{\texttt{a}}(h)$ is a unit of $\mK[h^{\pm 1}]$, then $M$ is isomorphic to $A/A(h-\gamma)$ for some $\gamma \in \mK^*$. 
		\item If $\emph{\texttt{a}}(h)=\mu h^i (h-\zeta)$ for some $\mu, \zeta \in \mK^*$ and $i\in \Z$, then
		$M$ is isomorphic to one of the following modules
		\begin{itemize}[noitemsep,label=\raisebox{0.35ex}{\tiny$\bullet$}]
			\item $A/A(h-\gamma)$ where $\gamma \in \mK^*$ and $\gamma \notin \mathscr{O}(\zeta)$;
			\item $A/A(h-\zeta, \,x)$;
			\item $A/A(h-q^{-1}\zeta, \,y)$.
		\end{itemize}
		In particular, $\GK(M)=1$. 
	\end{enumerate}	
\end{proposition}
\begin{proof}
	1. In this case the algebra $A$ is isomorphic to the quantum torus $\mK_q[x^{\pm 1}, h^{\pm 1}]$. Since $M$ is $\mK[h^{\pm 1}]$-torsion, there exists a non-zero element $m\in M$ such that $(h-\gamma)m=0$ for some $\gamma \in \mK^*$. Then there is a natural epimorphism of $A$-modules $\pi: A/A(h-\gamma) \twoheadrightarrow M$. It is easy to see that the module $A/A(h-\gamma)$ is simple, hence $\pi$ is an isomorphism. 
	
	2. First, we show that the modules described in the proposition are simple modules.  If we denote $W(\gamma):=A/A(h-\gamma)$ and $\bar{1}=1+A(h-\gamma)$, then using the $\Z$-grading of $A$ we obtain  $W(\gamma)=\mK[x] \bar{1}\oplus \mK[y]\bar{1}$. Since $\gamma \notin \mathscr{O}(\zeta)$, the simplicity of $W(\gamma)$ then follows from the following identities: for $k\geq 1$,
	\begin{subequations}
		\begin{align}
		xy^k \bar{1} &=\texttt{a}(q^k \gamma) y^{k-1} \bar{1} \in \mK^* y^{k-1} \bar{1},  \label{xyk1}\\ 
		yx^k \bar{1} &=\texttt{a}(q^{-k+1}\gamma)x^{k-1} \bar{1}\in \mK^* x^{k-1} \bar{1}.  \label{xyk2} 
		\end{align}
	\end{subequations}
	Note that $W:=A/A(h-\zeta,x)=\mK[y]\tilde{1}$ where $\tilde{1}=1+A(h-\zeta,x)$. The simplicity of $W$ then follows from the identity (\ref{xyk1}) with $\gamma=\zeta$. Similarly, $W':=A/A(h-q^{-1}\zeta,y)=\mK[x] \tilde{\bf{1}}$ where $\tilde{\bf{1}}=1+A(h-q^{-1}\zeta,y)$. The simplicity of $W'$ follows from the identity (\ref{xyk2}) with $\gamma=q^{-1}\zeta$. 
	
	Since $M$ is $\mK[h^{\pm 1}]$-torsion and the field $\mK$ is algebraically closed, there exists a non-zero element $m\in M$ such that $(h-\gamma)m=0$ for some $\gamma \in \mK^*$. In particular, $M$ is an epimorphic image of $W(\gamma)$. If $\gamma \notin \mathscr{O}(\zeta)$,  the simplicity of $W(\gamma)$ then yields that $M \simeq W(\gamma)$. 	
	Now assume that $\gamma \in \mathscr{O}(\zeta)$, that is, $\gamma=q^k \zeta$ for some $k\in \Z$. If $k \geq 0$ then $e:=x^k m$ is a non-zero element of $M$. There are two cases: either $xe =0$ and $h e=\zeta e$, or otherwise, $f:=xe=x^{k+1}m \neq 0$ and then $yf=0$ and $hf=q^{-1}\zeta f$. In the first case, $M$ is an epimorphic image of $W$, and then the simplicity of $W$ yields that $M \simeq W$. In the second case, $M$ is an epimorphic image of $W'$ and the simplicity of $W'$ implies that $M \simeq W'$. For $k<0$, a similar argument shows that $M$ is either isomorphic to $W$ or to $W'$. It is clear that the modules $W(\gamma)$, $W$ and $W'$ have GK dimension one. 
\end{proof}

We next turn to the description of simple $\mK[h^{\pm 1}]$-torsionfree $A$-modules. 
Let $\CB$ be the localization of the algebra $A=A(\texttt{a}(h),q)$ at the Ore set $S=\mK[h^{\pm 1}]\setminus\{0\}$. Then $\CB=AS^{-1}=\mK(h)[x^{\pm 1};\s]$ is a skew polynomial algebra where $\mK(h)$ is the field of rational functions in $h$ and $\s(h)=qh$.  We note that the localization of $A$ at the powers of the element $x$ is a quantum torus, that is, $A[x^{-1}]=\mK[h^{\pm 1}][x^{\pm 1};\s]$ where $\s(h)=qh$. Thus we have the following inclusions of algebras: $A \subset A[x^{-1}] \subset \CB$. The algebra $\CB$ is a Euclidean ring with left and right division algorithms. In particular, $\CB$ is a principal left and right ideal domain. By a classical theory, the simple $\CB$-modules can be described in terms of irreducible elements of $\CB$ (that is, if $f=p_1p_2 \in \CB$ then either $p_1$ or $p_2$ is a unit of $\CB$). More precisely, a $\CB$-module $M$ is simple iff $M \simeq \CB/\CB f$ for some irreducible element $f \in \CB$, see for instance \cite{Jacobson}.

An element $f=y^m \beta_m+\cdots+y\beta_1+\beta_0 \in A$ of length $m >0$ where $\beta_i \in \mK[h^{\pm 1}]$, is \emph{$\ell$-normal} if the following condition holds: every scalar $q^j \zeta$ ($j\geq 0$) is not a root of $\beta_0$, and if $\l$ and $q^{-k} \l$ are roots of $\beta_0$ and $\beta_m$ respectively then $k>0$. 

The following proposition gives a description of the set $\widehat{A}\,(\mK[h^{\pm 1}]\text{-torsionfree})$.
\begin{proposition} \label{a20Nov21} 
	(\cite{Bav-SimpleD[xy]Mod,Bav-VanOys-SimpModCrossProd}).
	Let $M$ be a simple $\mK[h^{\pm 1}]$-torsionfree module over the algebra $A=A(\emph{\texttt{a}}(h),q)$. Then $M \simeq A/A \cap \CB f$ for an irreducible $\ell$-normal element $f=y^m\beta_m+\cdots +\beta_0$ of $\CB$.	
\end{proposition}
In the following example, we give two families of simple $\mK[h^{\pm 1}]$-torsionfree $A$-modules. 
\begin{example} \label{b19Nov21} 
	For any $\gamma \in \mK^*$, we define 
	\begin{equation*}
	\mX(\gamma):=A/A(x-\gamma), \quad {\rm and}\quad \mathbb{Y}(\gamma):=A/A(y-\gamma). 
	\end{equation*}
	Then $\mX(\gamma)$ and $\mathbb{Y}(\gamma)$ are simple $\mK[h^{\pm 1}]$-torsionfree $A$-modules of GK dimension 2. 
\end{example}
\begin{proof}
	We only prove the result for $\mX(\gamma)$, since the proof for $\mathbb{Y}(\gamma)$ is similar. 
	Using the $\Z$-grading of $A$, we obtain $\mX(\gamma)=\oplus_{i\in \N} \mK[h^{\pm 1}]y^i \bar{1}$ where $\bar{1}=1+A(x-\gamma)$. From this it is clear that the module $\mX(\gamma)$ is $\mK[h^{\pm 1}]$-torsionfree with $\GK(\mX(\gamma))=2$. To show that $\mX(\gamma)$ is simple, let $V$ be a non-zero submodule of $\mX(\gamma)$. Note that for $i \geq 1$, $xy^i \bar{1}=\texttt{a}(qh)y^{i-1} \bar{1}$. Using this identity and the fact that $h$ is a normal element of $A$, we deduce that there is a non-zero element $v\in \mK[h^{\pm 1}]\bar{1}$ that belongs to $V$. Write $v=\sum_{i=m}^n \alpha_i h^i \bar{1}$ where $\alpha_i\in \mK$ and $m,n\in \Z$. Notice that the summands $h^i \bar{1}$ are $x$-eigenvectors with distinct eigenvalues ($x\cdot h^i \bar{1}=q^i \gamma h^i \bar{1}$). It follows that $h^i \bar{1} \in V$ if $\alpha_i \in \mK^*$. Thus $\bar{1}\in V$ and $V=\mX(\gamma)$.  This proves the simplicity of $\mX(\gamma)$. 
\end{proof}

\begin{remark}
	According to (\ref{DAmod}), Propositions \ref{a19Nov21} and    \ref{a20Nov21} give a classification of simple $A$-modules (up to irreducible elements of some Euclidean ring).  In view of (\ref{Cchi}) and (\ref{Cchip}), this gives a classification of simple unfaithful $\sC(\chi)$-modules. Applying Theorem \ref{D17Nov21}, we obtain a great number of simple weight $D_q$-modules. 
\end{remark}

There are plenty of simple faithful $\sC(\chi)$-modules, however,
it would be a very difficult task to give a classification of all such modules. We end this section with several examples of simple faithful $\sC(\chi)$-modules, and the description of the corresponding simple weight $D_q$-modules. 
\begin{proposition} \label{20Nov21} 
	For any $\chi \in \mK^*$, we define
	\begin{equation*}
	\begin{aligned}
	H(\chi):&= \sC/\sC(K-\chi, x_1,y_2), &L(\chi):&= \sC/\sC(K-\chi,y_1,x_2),\\
	M(\chi):&= \sC/\sC(K-\chi, x_1,x_2), &N(\chi):&= \sC/\sC(K-\chi,y_1,y_2). 
	\end{aligned}
	\end{equation*}
	\begin{enumerate}
		\item The modules $H(\chi)$, $L(\chi)$, $M(\chi)$ and $N(\chi)$ are simple faithful $\sC(\chi)$-modules. 
		\item The modules $\ind_{\sC}^{D_q}(H(\chi))$, $\ind_{\sC}^{D_q}(L(\chi))$, $\ind_{\sC}^{D_q}(M(\chi))$ and $\ind_{\sC}^{D_q}(N(\chi))$ are simple faithful weight $D_q$-modules. Moreover, 
		\begin{equation*}
		\begin{aligned}
		\ind_{\sC}^{D_q}(H(\chi)) &\simeq D_q/D_q(K-\chi, E,b), \quad  &\ind_{\sC}^{D_q}(L(\chi)) \simeq D_q/D_q(K-\chi, F,c), \\
		\ind_{\sC}^{D_q}(M(\chi)) &\simeq D_q/D_q(K-\chi, E, F), 	&\ind_{\sC}^{D_q}(N(\chi)) \simeq D_q/D_q(K-\chi, b,c). 
		\end{aligned}
		\end{equation*}
	\end{enumerate}
\end{proposition}
\begin{proof}
	We only prove the results for the modules $H(\chi)$ and $\ind_{\sC}^{D_q}(H(\chi))$, since other cases can be proved similarly.
	
	1. Note that the $\sC$-module $H(\chi)$ is naturally viewed as a $\sC(\chi)$-module, since $K$ is central in $\sC$.  Using the basis of $\sC$, we obtain $H(\chi)=\oplus_{i,j\in \N}\mK y_1^i x_2^j \bar{1}$	where $\bar{1}=1+\sC(K-\chi, x_1,y_2)$. The simplicity of $H(\chi)$ then follows from the following two equalities in $\sC(\chi)$, 
	\begin{equation*}
	y_2 x_2^j=q^{2j}x_2^j y_2-q^3 \tfrac{1-q^{2j}}{1-q^2}x_2^{j-1}, \quad x_1 y_1^i=q^{2i}y_1^i x_1+\chi \tfrac{1-q^{2i}}{1-q^2}y_1^{i-1}. 
	\end{equation*}
	In particular,  the annihilator of the module $H(\chi)$ in the algebra $\sC(\chi)$ is a primitive ideal of $\sC(\chi)$. Notice that $u\bar{1}=\chi \bar{1}$ and $v\bar{1}=q^2 \bar{1}$, thus the elements $u$ and $v$ do not belong to the annihilator of $H(\chi)$. Then from the classification of primitive ideals of $\sC(\chi)$, see Theorem \ref{A17Nov21}.(3), it follows that $\ann_{\sC(\chi)}H(\chi)=0$. This proves that $H(\chi)$ is a faithful $\sC(\chi)$-module. 
	
	2. The fact that $\ind_{\sC}^{D_q}(H(\chi))$ is a simple weight $D_q$-modules follows immediately from Theorem \ref{D17Nov21}. The isomorphisms in the statement are evident. If we denote by $\bar{1}=1+D_q(K-\chi,E,b)$ the canonical generator of the module $M=D_q/D_q(K-\chi,E,b)$, then an easy calculation shows that $\phi \bar{1}=q^2 a\bar{1}$ and $\psi \bar{1}=\chi a^{-1} \bar{1}$. Since $\ann_{D_q}(M)$ is a primitive ideal of $D_q$ and the elements $\phi$ and $\psi$ do not belong to $\ann_{D_q}(M)$, it follows from Corollary \ref{B19Oct21}.(3) that $\ann_{D_q}(M)=0$. It means that the module $M$ (and hence $\ind_{\sC}^{D_q}(H(\chi))$) is a  faithful $D_q$-module. 
\end{proof}

Note that the module $\ind_{\sC}^{D_q}(M(\chi))$ is isomorphic to $\CH(\chi):=D_q \otimes_{U_q}U_q/U_q(K-\chi,E,F)$, that is,  the module induced from a one-dimensional $U_q$-module, see Corollary \ref{d16Oct21}. As a vector space, $\CH(\chi)$ is isomorphic to $\CO_q$. In the case $\chi=1$, it is known as the Heisenberg representation.

\section{Simple $a$-weight modules of $D_q(E_2)$} \label{aWeightMod} 
In this section, we consider a class of simple $D_q$-modules which are $\mK[K,K^{-1}]$-torsionfree, and hence are not weight modules studied in the previous section.
Since the element $a$ is normal in $D_q$, similar to the previous section, one may consider simple $\mK[a, a^{-1}]$-torsion $D_q$-modules; that is, the simple $a$-weight modules defined below. 
\begin{definition}
	A $D_q$-module $M$ is called an \emph{$a$-weight module} provided that $M=\oplus_{\mu \in \mK^*}M_{\mu}$ where $M_{\mu}=\{ m \in M\,|\, am=\mu m \}$. We denote by $\supp_a(M):=\{ \mu\in \mK^*\,|\, M_{\mu} \neq 0 \}$ the set of all $a$-weights of $M$. 
\end{definition}

From the commutation relations of the element $a$ with the generators of $D_q$, for an $a$-weight module $M$ with $\mu \in \supp_a(M)$, the following hold
\begin{equation*}
\begin{aligned}
&a: M_{\mu} \rightarrow M_{\mu}, &&b: M_{\mu} \rightarrow M_{q\mu}, &&c: M_{\mu} \rightarrow M_{q\mu}, \\
&K: M_{\mu} \rightarrow M_{q\mu}, &&E: M_{\mu} \rightarrow M_{\mu}, &&F: M_{\mu} \rightarrow M_{q^{-1}\mu}. 
\end{aligned}
\end{equation*}

The following proposition describes the centralizer $C(a, D_q)$ of the element $a$ in the algebra $D_q$. In particular, it shows that $C(a, D_q)$ is isomorphic to the algebra $C(K, D_q)$ which has been described in Proposition \ref{a16Nov21}. Here we choose a particular set of generators of $C(a, D_q)$ to make this isomorphism evident. 
\begin{proposition} \label{a21Nov21} 
	Let $w:=q^{-1}a^2K^{-1}c$, $t_1:=KF$ and $t_2:=q^3 a^{-1}K^{-1}b$. 
	\begin{enumerate}
		\item The algebra $\mathscr{A}:=C(a, D_q)$ is generated by $a^{\pm 1}$, $E$, $w$, $t_1$, $t_2$ subject to the following defining relations	
		\begin{equation} \label{CenAD} 
		\begin{aligned}
		&Ew-q^2 wE=a, &\quad  &t_1 t_2-q^{-2}t_2 t_1=q, \\
		&t_1E=q^{2}Et_1, & &t_2E=q^{-2}Et_2, \\
		&t_1 w=q^{-2} w t_1, &&t_2w=q^{2} wt_2, \\
		& a a^{-1}=1=a^{-1}a, && [a,\,\, \cdot]=0. 
		\end{aligned}
		\end{equation}
		In particular, $\mathscr{A}$ is a Noetherian domain of \emph{GK} dimension 5. 
		\item The algebra $D_q$ is free as a left and right $\mathscr{A}$-module, $D_q=\bigoplus_{i\in \Z} \mathscr{A} K^i = \bigoplus_{i\in \Z} K^i \mathscr{A}$. 
		\item The following map defines an isomorphism of algebras 
		\begin{equation*}
		\Phi: \mathscr{A} \rightarrow \sC,  \quad a^{\pm 1} \mapsto K^{\pm 1},\;\; E \mapsto x_1, \;\;w\mapsto y_1, \;\;t_1 \mapsto x_2, \;\;t_2 \mapsto y_2. 
		\end{equation*}
	\end{enumerate}	
\end{proposition}
\begin{proof}
	1. The adjoint action of the element $a$ on $D_q$ gives rise to a $\Z$-grading of $D_q$, and the algebra $\mathscr{A}$ consists of homogeneous elements of degree zero with respect to this grading. Thus
	\begin{equation*}
	\begin{aligned}
	\mathscr{A}=&\bigoplus_{\stackrel{i+k_2-k_3+k_4=0}{i,j\in\Z, k_1,\ldots,k_4\in\N}} \mK K^i a^j E^{k_1}c^{k_2}F^{k_3}b^{k_4} \\
	=&\,\, \bigoplus_{j\in \Z, k_1, \ldots, k_4\in \N} \mK K^{-k_2+k_3-k_4}a^j E^{k_1} c^{k_2} F^{k_3}b^{k_4}\\
	=&\,\, \bigoplus_{j\in \Z, k_1, \ldots, k_4\in\N} \mK a^j E^{k_1}(K^{-1}c)^{k_2}(KF)^{k_3}(K^{-1}b)^{k_4}\\
	=&\,\, \bigoplus_{j\in\Z, k_1,\ldots,k_4\in\N} \mK a^j E^{k_1}w^{k_2}t_1^{k_3}t_2^{k_4}. 		
	\end{aligned}
	\end{equation*}
	Thus $\mathscr{A}$ is generated by the elements $a^{\pm 1}$, $E$, $w$, $t_1$, $t_2$	 and the elements $a^j E^{k_1}w^{k_2}t_1^{k_3}t_2^{k_4}$ (where $j\in\Z$, $k_1,\ldots,k_4\in\N$) form a basis of $\mathscr{A}$. It is straightforward to verify that the relations in (\ref{CenAD}) hold in the algebra $\mathscr{A}$. Let $R$ be the abstract $\mK$-algebra generated by $a^{\pm 1}$, $E$, $w$, $t_1$ and $t_2$ subject to the relations in (\ref{CenAD}). Then there is a natural epimorphism of algebras $\pi:R \rightarrow \mathscr{A}$. Comparing the bases of $R$ and $\mathscr{A}$ yields that $\pi$ is an isomorphism. It is easy to see that $\mathscr{A}$ can be presented as an iterated Ore extension, in particular, $\mathscr{A}$ is a Noetherian domain of GK dimension 5. 
	
	2. Statement 2 follows from the bases of $\mathscr{A}$ and $D_q$ given in the proof of statement 1. 
	
	3. By the defining relations of $\mathscr{A}$ and $\sC$ (see Proposition \ref{a16Nov21}), statement 3 is evident.	
\end{proof}

For any $\chi \in \mK^*$, define $\mathscr{A}(\chi):=\mathscr{A}/\mathscr{A}(a-\chi)$. Clearly, we have $\mathscr{A}(\chi) \simeq \sC(\chi)$. Using this isomorphism, the results obtained for the algebra $\sC$ (resp. $\sC(\chi)$) can be translated to that of $\mathscr{A}$ (resp. $\mathscr{A}(\chi)$). In particular, the centre of $\mathscr{A}$ is the Laurent polynomial algebra $\mK[a, a^{-1}]$ (see Lemma \ref{b16Nov21}), the prime and primitive spectra of $\mathscr{A}(\chi)$ are obtained (see Theorem \ref{A17Nov21}), the algebra $\mathscr{A}$ has no finite-dimensional modules (see Corollary \ref{c17Nov21}), and a classification of all simple unfaithful $\mathscr{A}(\chi)$-modules is obtained. 
The adjoint action of the normal element $K$ on $\mathscr{A}$ gives rise to an automorphism $\tau_K$ of $\mathscr{A}$; that is, $\tau_K: \mathscr{A} \rightarrow \mathscr{A}$, $t \mapsto K^{-1}tK$. On the generators of $\mathscr{A}$, $\tau_K$ takes the following form
\begin{equation} \label{AutK} 
\tau_K: \quad a \mapsto qa, \,\,\, E \mapsto q^{-2}E, \,\,\, w \mapsto q^3 w, \,\,\, t_1 \mapsto q^2 t_1, \,\,\,t_2 \mapsto q^{-2}t_2. 
\end{equation}

The following theorem shows that the problem of classification of simple $a$-weight $D_q$-modules is essentially reduced to the problem of classification of all simple $\mathscr{A}$-modules. 
\begin{theorem} \label{A22Nov21} 
	\
	\begin{enumerate}
		\item 	Let $V$ be a simple $a$-weight $D_q$-module. Then $\supp_a(V)=\{ q^i \mu\,|\, i\in \Z \}$ for any $\mu \in \supp_a(V)$, and every $a$-weight space of $V$ is infinite-dimensional. 
		\item 	The following map defines a surjection
		\begin{equation*}
		\widehat{\mathscr{A}} \twoheadrightarrow \widehat{D_q}(a\text{-}{\rm weight}), \quad M \mapsto \ind_{\mathscr{A}}^{D_q}(M):=D_q \otimes_{\mathscr{A}} M. 
		\end{equation*}
		Moreover, for $M, M' \in \widehat{\mathscr{A}}$, $\ind_{\mathscr{A}}^{D_q}(M) \simeq \ind_{\mathscr{A}}^{D_q}(M')$ iff $_{\mathscr{A}} M' \simeq {_{\mathscr{A}}} M^{\tau_K^i}$ for some $i\in \Z$ where $\tau_K$ is the automorphism of $\mathscr{A}$ defined in (\ref{AutK}). 
	\end{enumerate}
\end{theorem}
\begin{proof}	
	1. The proof is similar to that of Theorem \ref{B17Nov21}. 
	
	2. The result can be proved in much the same way as Theorem \ref{D17Nov21}. 
\end{proof}

\noindent\textbf{Acknowledgements.}
This research was partially supported by the NSFC grant 11601167.

\small{ }

	\

School of Mathematical Science, Yangzhou University, Yangzhou, 225002, China

Email: wqtao@yzu.edu.cn

\end{document}